\documentclass[11pt, a4paper]{amsart}
\usepackage[latin1]{inputenc}
\usepackage{amsthm}
\usepackage{amssymb}
\usepackage{amsmath}
\usepackage{cancel}
\usepackage{a4wide}

\usepackage{lipsum}

\usepackage{verbatim}
\usepackage[normalem]{ulem}
\numberwithin{equation}{section} 

\theoremstyle{plain}
\newtheorem{thm}{Theorem}[section]

\theoremstyle{definition}
\newtheorem{defn}[thm]{Definition}
\theoremstyle{plain}
\newtheorem{prop}[thm]{Proposition}

\newtheorem*{notation}{Notation}

\newtheorem{lem}[thm]{Lemma}
\newtheorem{cor}[thm]{Corollary}
\theoremstyle{remark}
\newtheorem{rem}[thm]{Remark}
\newtheorem{ex}[thm]{Example}

\newtheorem{Hypothesis}[thm]{Hypothesis}

\usepackage{graphicx}
\usepackage{amscd}
\usepackage[all,knot,poly]{xy}
\usepackage[english]{babel}
\usepackage{xcolor}

\usepackage{mathrsfs}
\usepackage{hyperref}
\newcommand{\F}{\mathbb{F}}

\newcommand{\C}{\mathbb{C}}

\newcommand{\Z}{\mathbb{Z}}
\newcommand{\N}{\mathbb{N}}

\newcommand{\ones}{J}
\newcommand{\sJ}{\mathcal{J}}
\newcommand{\Rad}{{\rm Rad}}
\newcommand{\se}{\text{ss}}

\title{On the joins of group rings}
\author[S.K. Chebolu, J. Merzel, J. Min\'a\v{c}, L. Muller, T.T. Nguyen, F.W.~Pasini, N.D. Tan]{Sunil K. Chebolu, Jonathan L. Merzel, J\'an Min\'a\v{c}, Lyle Muller, \\ Tung T. Nguyen, Federico W. Pasini, Nguy$\tilde{\text{\^{e}}}$n Duy T\^{a}n}

\address{
Illinois State University}
\email{schebol@ilstu.edu}

\address{Soka University of America}
\email{jmerzel@soka.edu }

\address{
University of Western Ontario}
\email{minac@uwo.ca}

\address{University of Western Ontario}
\email{lmuller2@uwo.ca}

\address{University of Western Ontario and Onepick Inc}
\email{tungnt@uchicago.edu}

\address{Huron University College}
\email{fpasini@uwo.ca }

\address{Hanoi University of Science and Technology}
\email{tan.nguyenduy@hust.edu.vn}

\date{}

\begin{document}
\thanks{Sunil Chebolu is partially supported by Simons Foundation's Collaboration Grant for Mathematicians (516354). J\'an Min\'a\v{c} is partially supported by the Natural Sciences and Engineering Research Council of Canada (NSERC) grant R0370A01. J\'an Min\'a\v{c} also gratefully acknowledges 
Faculty of Sciences Distinguished Research Professorship award for 2020/21. J\'an Min\'a\v{c}, Lyle Muller, Tung T Nguyen, and Federico Pasini acknowledge the support of the Western Academy for Advanced Research. Nguy$\tilde{\text{\^{e}}}$n Duy T\^{a}n is funded by Vingroup Joint Stock Company and supported by Vingroup Innovation Foundation (VinIF) under the project code VINIF.2021.DA00030}


\dedicatory{Dedicated to Professor T. Y. Lam with gratitude and admiration.}

\maketitle

\begin{abstract}
Given a  collection $\{ G_i\}_{i=1}^d$ of finite groups and a ring $R$, we define a subring of the ring $M_n(R)$ ($n = \sum_{i=1}^d|G_i|)$  that encompasses all the individual group rings $R[G_i]$ along the diagonal blocks as $G_i$-circulant matrices. The precise definition of this ring was inspired by a construction in graph theory known as the joined union of graphs. We call this ring the join of group rings and denote it by $\mathcal{J}_{G_1,\dots, G_d}(R)$. In this paper, we present a systematic study of the algebraic structure of  $\mathcal{J}_{G_1,\dots, G_d}(R)$. We show that it has a ring structure and characterize its center,  group of units, and Jacobson radical. When $R=k$ is an algebraically closed field, we derive a formula for the number of irreducible modules over $\mathcal{J}_{G_1,\dots, G_d}(k)$.   We also show how a blockwise extension of the Fourier transform provides both a generalization of the Circulant Diagonalization Theorem to joins of circulant matrices and an explicit isomorphism between the join algebra and its Wedderburn components.

\end{abstract}

\tableofcontents
\section{Introduction}

Circulant matrices appear naturally in many problems in physics, spectral graph theory, and non-linear dynamics (see, for example, \cite{cir1, cir2, cir3, ko1}). They are closely connected with the theory of Fourier analysis and representation theory of finite groups. For example, for circulant matrices associated with a cyclic group, we can describe their spectrum explicitly via the discrete Fourier transform (see \cite{davis2013circulant} for an extensive treatment of this topic). For this reason, many problems involving circulant matrices have closed-form solutions.

Many real-world networks have structure beyond that of circulant networks. Let us imagine that there are $d$ networks with their own connections. These individual networks are not isolated. They interact via a modeled network $G$ in the following way. Suppose $G$ is a weighted  digraph with $d$ vertices $\{v_1, v_2, \ldots, v_d\}$.   Let $G_1, \ldots, G_d$ be weighted digraphs on pairwise disjoint sets of $k_1,  \ldots, k_d$ vertices. The joined union $G[G_1, \ldots, G_d]$ is obtained from the union of $G_1, \ldots, G_d$ by joining with an edge each pair of a vertex from $G_i$ and a vertex from $G_j$ whenever $v_i$ and $v_j$ are adjacent in $G$ (see \cite{joined_union} for further details). Let $A_{G} =(a_{ij})$ be the adjacency matrix of $G$ and $A_{G_1}, \ldots, A_{G_d}$ be the adjacency matrices of $G_1, \ldots, G_d$ respectively. We can then observe that the adjacency matrix of $G[G_1, \ldots, G_d]$ has the following form, which we refer to as a \textit{join of the matrices} $A_{G_1},\dots,A_{G_d}$
\begin{equation} \label{eq:joined_union1}
A=\left[\begin{array}{c|c|c|c}
A_{G_1} & a_{12} J_{k_1, k_2} & \cdots & a_{1d} J_{k_1, k_d} \\
\hline
a_{21} J_{k_2, k_1} & A_{G_2} & \cdots & a_{2d} J_{k_2, k_d} \\
\hline
\vdots & \vdots & \ddots & \vdots \\
\hline
a_{d1} J_{k_d, k_1} & a_{d2} J_{k_d, k_2} & \cdots & A_{G_d}
\end{array}\right].
\end{equation} 
Here $J_{m,n}$ is the matrix of size $m \times n$ with all entries equal to $1.$

In our investigation of the Kuramoto model of coupled oscillator networks, these joined networks appear quite often and provide some new interesting phenomena (\cite{ko3, CM1, CM1_b,ko4}). In particular, in \cite{CM1}, we describe the spectrum of the joins of circulant matrices (for cyclic groups) explicitly. We also apply our main results to study several edge-removal problems in spectral graph theory as well as describe new equilibrium points on the Kuramoto models (see \cite[Section 4, Section 5]{CM1}.)

For these reasons, it is important to develop a systematic understanding of the joins of circulant matrices, or more generally of matrices which are circulant with respect to a group $G$, according to Definition \ref{defn:G-circulant}. A crucial observation here is that, for fixed groups $G_1,\dots,G_d$, the set of all $A$ as in Equation \eqref{eq:joined_union1}, such that each $A_{G_i}$ is circulant with respect to $G_i$, has the structure of a ring with identity (we call it the \textit{join ring}). We can then utilize tools from ring theory and representation theory to study the abstract structural properties of the sets of such matrices. In this article, we lay the foundation for this line of research through a systematic study of the algebraic structure of this ring. Since every group algebra is also a join ring, the results of this paper can be viewed as natural generalizations of the corresponding results for group algebras. We now summarize our main results. We refer the reader to the main text for more precise statements.

\begin{thm} Let $R$ be a unital ring and let  $\mathcal{J}_{G_1,\dots, G_d}(R)$ denote the set of all matrices of the form \ref{eq:joined_union1}  where each $A_i$ is a $G_i$-circulant matrix over $R$. Then we have the following.
\begin{enumerate}
    \item $\mathcal{J}_{G_1,\dots, G_d}(R)$ has the structure of a unital ring and there is an augmentation map $\epsilon \colon \mathcal{J}_{G_1,\dots, G_d}(R) \rightarrow M_d(R)$ that generalizes the augmentation map on group rings.
    \item If $R$ is a semisimple ring and  $|G_i|$ is invertible in $R$ for all $1 \leq i \leq d$,  then the ring $\sJ_{G_1, G_2, \ldots, G_d}(R)$ is semisimple.
    
    \item If $R$ is a commutative ring such that $|G_i|=0$ in R for all $1 \leq i \leq d$, then
\[ (\mathcal{J}_{G_1,\dots, G_d}(R))^{\times} \cong (R^{d^2-d} \times \prod_{i=1}^d U_1(R[G_i])) \rtimes (R^{\times})^d ,\]
where $U_1(R[G_i])$ is the group of principal units in $R[G_i]$.
\item For $A$ as above (\ref{eq:joined_union1}) and $k$ any field, $A \in \Rad(\mathcal{J}_{G_1, G_2, \ldots, G_d}(k))$ if
and only if for all $i$ we have $A_{i}\in \Rad(k[G_{i}])$  and   
whenever $p\nmid |G_i||G_{j}|, 1\leq i\neq j\leq d$, we have $a_{ij}=0$. In particular,  $\sJ_{G_1, G_2, \ldots, G_d}(k)$ is semisimple if and only if $|G_i|$ are invertible in $k$. 
\item If $k$ is algebraically closed and $\text{char}(k)$ is relatively prime to $\prod_{i=1}^d |G_i|$, then the number of irreducible modules over $\sJ_{G_1, G_2, \ldots, G_d}(k)$ is 
$c(G_1)+c(G_2)+\ldots+c(G_d)-d+1,$ 
where $c(G_i)$ is the number of conjugacy classes of $G_i$.
\item If $k$ is any field and $d \geq 2$, then $\sJ_{G_1, G_2, \ldots, G_d}(k)$ is a Frobenius algebra if and only if $|G_i|$ is invertible in $k$ for all $1 \leq i \leq d.$ 

\end{enumerate}
\end{thm}

The structure of this article is as follows. In Section $2$, we recall the definition of a $G$-circulant matrix over a ring $R.$ We also provide several characterizations for $G$-circulant matrices. Using these characterizations, we reprove several results of Hurley (see \cite{hurley2006}) on the structure of the ring of $G$-circulant matrices. In Section 3, we introduce the ring of the joins of $G$-circulant matrices. We then discuss some of its basic properties, such as the existence of a generalized augmentation map, its center, and its decomposition in the semisimple case. In particular, we prove a generalized version of Maschke's theorem. Section 4 studies the unit group of the join algebra $\sJ_{G_1, G_2, \ldots, G_d}(k)$ over a field $k$. Here, we provide a complete characterization for a join matrix to be invertible in $\sJ_{G_1, G_2, \ldots, G_d}(k)$. Additionally, we study the precise structure of the unit group in the modular case (where $\text{char}(k)$ divides $|G_i|$ for all $i$). In Section 5, we determine the Jacobson radical of the join ring, as well as its semisimplification. In section 6, we discuss the necessary conditions for the join algebra  to be a Frobenius algebra. Finally, in the last section, we study the explicit structure of the join algebra $\sJ_{G_1, G_2, \ldots, G_d}(k)$ when $k$ is a field and the groups $G_i$ are all cyclic.

\section*{Acknowledgements}
We thank professor T. Y. Lam for several comments and suggestions that helped improve the exposition of this paper.

\section{$G$-circulant matrices over an arbitrary ring}
Let $R$ be a ring (with unity), fixed for this section. \ A circulant $%
n\times n$ matrix over $R$ is a matrix of the form 
\[
\left[ 
\begin{array}{lllll}
a_{1} & a_{2} & \cdots  & a_{n-1} & a_{n} \\ 
a_{n} & a_{1} & \cdots  & a_{n-2} & a_{n-1} \\ 
\vdots  & \vdots  & \ddots  & \vdots  & \vdots  \\ 
a_{3} & a_{4} & \cdots  & a_{1} & a_{2} \\ 
a_{2} & a_{3} & \cdots  & a_{n} & a_{1}%
\end{array}%
\right] 
\]%
where $a_{1},\ldots ,a_{n}\in R$. \ As we discussed in the introduction, such matrices have a wide variety of
applications in mathematics.  It is not difficult to show that the set of all circulant $n\times n$ matrices over $R$ forms a ring isomorphic to the group ring $R[G]$ where $G= \langle g \rangle$ is cyclic of order $n$ and the isomorphism takes $\sum\limits_{i=1}^{n}a_{i}g^{i-1}$ to the matrix above.

More generally, let $G$ be any finite group, say of order $n$, and fix an indexed listing of $G$ so that $G=\{g_{i}\}_{i=1}^{n}$.  It will be convenient for the purposes of this section to view $n\times n$ matrices over $R$ as having their rows and columns indexed by the elements of $G$ (so that the $i,j$ entry $M_{ij}$ of a matrix $M$ will be renamed the $%
g_{i},g_{j}$ entry).  

\begin{defn}\label{defn:G-circulant}
An $n \times n$ $G$-circulant matrix over $R$ is an $n\times n$ matrix 
\[
\left[ 
\begin{array}{llll}
a_{g_{1},g_{1}} & a_{g_{1},g_{2}} & \cdots  & a_{g_{1},g_{n}} \\ 
a_{g_{2},g_{1}} & a_{g_{2},g_{2}} & \cdots  & a_{g_{2},g_{n}} \\ 
\vdots  & \vdots  & \ddots  & \vdots  \\ 
a_{g_{n},g_{1}} & a_{g_{n},g_{2}} & \cdots  & a_{g_{n},g_{n}}%
\end{array}%
\right] 
\]%
over $R$ with the property that
for all $g,g_i,g_j \in G,$ $a_{g_i,g_j}=a_{gg_i,gg_j}$.\ 
\end{defn}

This is a (harmless) abuse of language, since the choice of an ordering on the elements of $G$ affects which matrices are called $G$-circulant. Whenever we refer to $G$-circulant matrices, we will assume that an ordering has been chosen once and for all on $G$.  Note that if $G=\{g_{\sigma(i)}\}_{i=1}^{n}$ is a reordering of the elements of $G$, then a matrix $A$ is $G$-circulant with respect to $\{g_{i}\}_{i=1}^{n}$ if and only if $PAP^{-1}$ is $G$-circulant with respect to $\{g_{\sigma(i)}\}_{i=1}^{n}$ where $P$ is the permutation matrix $(\delta_{i,\sigma(j)})$ and $\delta$ is the Kronecker delta.
\begin{rem}\label{rem:semimagic} It is immediate from the definition 
(using $a_{g_1,g_k}=a_{g_i,g_ig_1^{-1}g_k}=a_{g_jg_k^{-1}g_1,g_j}$) that the entries in any row or column of a $G$-circulant matrix are permutations of the elements of the first row. In particular, any $G$-circulant matrix is a semimagic square; i.e., all rows and columns have the same sum. We refer to \cite{murase1957semimagic} for further discussions about the ring of semimagic squares.

\end{rem}

A circulant matrix as above is then the special case where  $G= \langle g \rangle$ is cyclic of order $n$ (relative to choosing $g_{i}=g^{i-1})$.  $G$-circulant
matrices are defined and studied in \cite{hurley2006} and generically in \cite{kanemitsu2013matrices}; we will give alternate proofs for some of the theorems of those papers. In the process, some new results will emerge.

\begin{notation}
For $g\in G$ we denote by $P_{g}$ the permutation matrix such that right
multiplication by $P_{g}$ permutes columns (indexed by $g_{1},g_{2},\ldots $%
) by left multiplication by $g$.  That is to say, for $A\in M_{n,n}(R)$ the 
$gg_{j}$th column of $AP_{g\text{ }}$is the $g_{j}$th column of $A$. 
Explicitly, $(P_{g})_{g_{i},g_{j}}=\delta _{gg_{i},g_{j}}$(Kronecker delta).
\ We denote by $P_{g}^{\prime }$  the permutation matrix such that right
multiplication by $P_{g}^{\prime }$ permutes columns by \textit{right
multiplication by }$g$ (so the $g_{j}g$th column of $AP_{g}^{\prime }$ is
the $g_{j}$th column of $A$).  Explicitly $(P_{g}^{\prime
})_{g_{i},g_{j}}=\delta _{g_{i}g,g_{j}}$.  
\end{notation}

\begin{rem}
Notice that $g\neq h$ implies that the sets of positions where $P_{g}$ and $%
P_{h}$ have entries of $1$ are disjoint, and in particular then that $%
P_{g}\neq P_{h}$; the same considerations hold for $P_{g}^{\prime }$ and $%
P_{h}^{\prime }$.
\end{rem}

\begin{rem}
Using $g^{op}$ to work in the opposite group $G^{op}~$(and listing $%
G^{op}=\{g_{i}^{op}\}_{i=1}^{n}$), $P_{g}^{\prime }=P_{g^{op}}$.
\end{rem}

\begin{prop}
For all $g,h\in G$ we have $P_{g}P_{h}=P_{hg}$ and $P_{g}^{\prime
}P_{h}^{\prime }=P_{gh}^{\prime }$. 
\end{prop}

\begin{proof}
$(P_{g}P_{h})_{g_{r},g_{s}}=\sum\limits_{k=1}^{n}\delta
_{gg_{r},g_{k}}\delta _{hg_{k},g_{s}}=\delta _{hgg_{r},g_{s}}=\left(
P_{hg}\right) _{g_{r},g_{s}}$.  A similar computation establishes the
second equation, or we can apply the remark above.
\end{proof}

\begin{cor}
The maps $g\mapsto P_{g^{-1}}$ and $g\mapsto P_{g}^{\prime }$ give
isomorphisms from $G$ to $\{P_{g}\mid g\in G\}$ and to $\{P_{g}^{\prime }\mid g\in
G\}$ respectively.
\end{cor}

\begin{proof}
{}Immediate from the last proposition and our earlier remark implying that
these maps are one-to-one.
\end{proof}

\begin{prop}
\label{prop:characterization_G_circulant}
For any ring $R$, finite group $G=\{g_{i}\}_{i=1}^{n}$, and $A\in M_{n,n}(R)$
the following are equivalent.
\begin{enumerate}
\item $A$ is $G$-circulant.
\item  For all $g_i,g_j,g_k,g_l \in G$,  $g_{i}^{-1}g_{j}=g_{k}^{-1}g_{l}$ implies
that $A_{g_{i},g_{j}}=A_{g_{k},g_{l}}.$
\item There is some $\sum\limits_{k=1}^{n}b_{g_{k}}g_{k}\in R[G]$ such that 
$A_{g_{i},g_{j}}=b_{g_{i}^{-1}g_{j}}$ for all $1\leq i,j\leq n$.
\item There exist $c_{1},\ldots ,c_{n}\in R$ such that $A=\sum%
\limits_{i=1}^{n}c_{i}P_{g_{i}}^{\prime }$.
\item $P_{g}A=AP_{g}$ for all $g\in G$.
\end{enumerate}
\end{prop}

\begin{proof}
(1)$\Leftrightarrow $(2): $\ $If (1) holds and $g_{i}^{-1}g_{j}=g_{k}^{-1}g_{l}$ then $%
A_{g_{i},g_{j}}=A_{1,g_{i}^{-1}g_{j}}=A_{1,g_{k}^{-1}g_{l}}=A_{g_{k},g_{l}%
\text{.}}$ Conversely, if (2) holds, just note that $%
(gg_{i})^{-1}(gg_{j})=g_{i}^{-1}g_{j}$ to see that (1) holds as well. 

\medskip
\noindent (2)$\Rightarrow $(3):  For $1\leq k\leq n$ set $b_{k}=A_{g_{1},g_{1}g_{k}}$%
. \ Then $A_{g_{i},g_{j}}=A_{g_{1},g_{1}g_{i}^{-1}g_{j}}=b_{g_{i}^{-1}g_{j}}$%
.\newline

\smallskip
\noindent (3)$\Rightarrow $(4):  Letting $e_{g,h}$ denote the ``matrix unit" having entry $1$ in the $g,h$ position and $0$ elsewhere, we may write $A=\sum%
\limits_{g_{i},g_{j}}b_{g_{i}^{-1}g_{j}}e_{g_{i},g_{j}}$ . \ Setting $%
g_{k}=g_{i}^{-1}g_{j}$ and reindexing the sum, we have $A=\sum%
\limits_{g_{k},g_{j}}b_{g_{k}}e_{g_{j}g_{k}^{-1},g_{j}}=\sum%
\limits_{g_{k}}b_{g_{k}}\left(
\sum\limits_{g_{j}}e_{g_{j}g_{k}^{-1},g_{j}}\right) $.  We claim that $%
\sum\limits_{g_{j}}e_{g_{j}g_{k}^{-1},g_{j}}$ is just $P_{g_{k}}^{\prime }$, which gives the desired result.  To see this, we compute 
\[
\left( \sum\limits_{g_{j}}e_{g_{j}g_{k}^{-1},g_{j}}\right)
_{g_{r},g_{s}}=\sum\limits_{g_{j}}\delta _{g_{j}g_{k}^{-1},g_{r}}\delta
_{g_{j},g_{s}}=\delta _{g_{s}g_{k}^{-1},g_{r}}=\delta
_{g_{r}g_{k},g_{s}}=(P_{g_{k}}^{\prime })_{g_{r},g_{s}}\text{.}
\]%

\smallskip
\noindent (4)$\Rightarrow $(5):  By (4), it is sufficient to show that for all $%
g_{i},g_{j}$ we have $P_{g_{i}}^{\prime }$ commutes with $P_{g_{j}}$.  This
is essentially the fact that left multiplications in a group commute with
right multiplications, which is just associativity. The explicit
calculation:%
\[
\left( P_{g_{i}}^{\prime }P_{g_{j}}\right)
_{g_{r},g_{s}}=\sum\limits_{g_{t}}\delta _{g_{r}g_{i},g_{t}}\delta
_{g_{j}g_{t},g_{s}}=\delta _{g_{j}(g_{r}g_{i}),g_{s}}=\delta
_{(g_{j}g_{r})g_{i},g_{s}}=\sum\limits_{g_{t}}\delta
_{g_{j}g_{r},g_{t}}\delta _{g_{t}g_{i},g_{s}}=\left(
P_{g_{j}}P_{g_{i}}^{\prime }\right) _{g_{r},g_{s}}.
\]%

\noindent (5)$\Rightarrow $(1):  We have 
\begin{eqnarray*}
\left( P_{g}A\right) _{g_{r},g_{s}} &=&\sum\limits_{g_{k}\in
G}(P_{g})_{g_{r},g_{k}}A_{g_{k},g_{s}}=\sum\limits_{g_{k}\in G}\delta
_{gg_{r},g_{k}}A_{g_{k},g_{s}}=A_{gg_{r}.g_{s}}\text{ \ and } \\
\left( AP_{g}\right) _{g_{r},g_{s}} &=&\sum\limits_{g_{k}\in
G}A_{g_{r},g_{k}}(P_{g})_{g_{k},g_{s}}=\sum\limits_{g_{k}\in
G}A_{g_{r},g_{k}}\delta _{gg_{k},g_{s}}=A_{g_{r},g^{-1}g_{s}}
\end{eqnarray*}%
\newline
so if $A$ and $P_{g}$ commute for all $g$ then for all $g_{r},g_{s},g$ we
have  $A_{g_{r},g^{-1}g_{s}}=A_{gg_{r}.g_{s}}$. \ Momentarily fixing $g$ and 
$g_{r}$ and setting $g_{t}=g^{-1}g_{s}$, $g_{t}$ runs through $G$ as $g_{s}$
does, and we then have $A_{g_{r},g_{t}}=A_{gg_{r}.gg_{t}}$, establishing (1).
\end{proof}

\begin{cor}[\cite{hurley2006}] \label{cor:groupring} 
 The set of all $G$-circulant matrices over $R$ forms a ring
isomorphic to the group ring $R[G]$.
\end{cor}

\begin{proof}
From our earlier observation that the permutation matrices $P_{g}^{\prime }$
are $\{0,1\}$-matrices, no two of which have entries of $1$ in the same
position, combined with the equivalence of condition (1) and condition (4)
in the previous proposition, we have that the set of all $G$-circulant
matrices over $R$ is a free left $R$-module on generators $\{P_{g}^{\prime
}\mid g\in G\}$; by our prior proposition that $P_{g}^{\prime }P_{h}^{\prime
}=P_{gh}^{\prime }$ it now follows that the map $\sum%
\limits_{k=1}^{n}b_{k}g_{k}\mapsto \sum\limits_{k=1}^{n}b_{k}P_{g\text{ }%
}^{\prime }$is the desired isomorphism.
\end{proof}

\begin{cor}[\cite{hurley2006}] \label{Hurley} \ If a $G$-circulant matrix $A$ is invertible as an element of $%
M_{n,n}(R)$, then it is a unit in the ring of $G$-circulant matrices.
\end{cor}

\begin{proof}
Since $A$ commutes with all $P_{g}$, so does its matrix inverse $A^{-1}$.
\end{proof}

Let $Z(R)$ denote the center of the ring $R$.

\begin{cor} The centralizer of the ring of $G$-circulant matrices (with respect to the listing $G=\{g_{i}\}_{i=1}^{n}$) in $M_{n,n}(R)$ is exactly the ring of $G^{op}$-circulant matrices (with respect to the listing $G^{op}=\{g_{i}^{op}\}_{i=1}^{n}$) in $M_{n,n}(Z(R))$. 
\end{cor}
\begin{proof} This follows from the equivalence of (1), (4), and (5) in the above
proposition and our earlier observation that $P_{g}^{\prime }=P_{g^{op}}$.
\end{proof}

\begin{lem}
\label{lem:transpose_circulant} If $A$  is  $G$-circulant then its transpose $A^T$ is also $G$-circulant.  
\end{lem}
\begin{proof} Suppose that $A$ is $G$-circulant. For every $g,g_i,g_j\in G$, one has
\[
 (A^T)_{gg_i,gg_j}=A_{gg_j,gg_i}=A_{g_j,g_i}=(A^T)_{g_i,g_j}. 
\]
Hence $A^T$ is $G$-circulant.
\end{proof}

\section{The ring of joins of $G$-circulant matrices}

\begin{defn}
Let $R$ be a (unital, associative) ring, $G_{1},\ldots ,G_{d}$  finite
groups of respective orders $k_{1},\ldots ,k_{d}$, and let $C_{i}$ be $G_{i}$%
-circulant ($1\leq i\leq d$) over $R$. \ By a join of $C_{1},\ldots ,C_{d}$
over $R$, we mean a matrix of the form 

\begin{equation}
\label{eq:join circulant matrix}
A=\begin{bmatrix} 
C_1&a_{12} J_{k_1,k_2} &\cdots & a_{1d}J_{k_1,k_d}\\
a_{21}J_{k_2,k_1} &C_2 &\cdots & a_{2d}J_{k_2,k_d}\\
\vdots&\vdots& &\vdots\\
a_{d1}J_{k_d,k_1}&a_{d2}J_{k_d,k_2}&\cdots& C_d
\end{bmatrix}
, 
\end{equation}
where $a_{ij}\in R\ (1\leq i\neq j\leq d)$ and $J_{r,s}$ denotes the $%
r\times s$ matrix, all of whose entries are $1\in R$. \ 

The join group ring of $G_{1},\ldots ,G_{d}$ over $R$, denoted $\sJ_{G_{1},\ldots
,G_{d}}(R)$, is the set of all such joins as the $C_{i\text{ }}$vary
independently through all $G_{i}$-circulant matrices ($1\leq i\leq d$) and
the $a_{ij}$ vary independently through all elements of $R$ ($1\leq i\neq
j\leq d$). \ (That this is a ring will be shown momentarily.)
\end{defn}

Let $M_n(R)$ denote the algebra of $n\times n$ matrices over $R$ where $n=\sum_{i=1}^d k_i$. We first have the following observation. Here we simply write $J$ for a matrix all of whose entries are $1\in R$ and whose dimensions can be inferred from the context (e.g., from a block structure).

\begin{prop}
For any ring $R$, finite groups $G_1,\ldots,G_d$, and $A\in M_{n,n}(R)$
the following are equivalent.
\begin{enumerate}
\item $A$ is in $\sJ_{G_1,\ldots,G_d}(R)$.
\item ${\rm diag}(P_{g_1},\ldots,P_{g_d})\cdot A=A\cdot {\rm diag}(P_{g_1},\ldots,P_{g_d})$ for every $g_1\in G_1,\ldots,g_d\in G_d$. Here for $g_i\in G_i$, $i=1,\ldots,d$, we denote ${\rm diag}(P_{g_1},\ldots,P_{g_d})$  the diagonal block matrix with diagonal block entries $P_{g_i}$.
\end{enumerate}
\end{prop}

\begin{proof} We write
\[
A=\left[\begin{array}{c|c|c|c}
A_{11} & A_{12} & \cdots & A_{1d} \\
\hline
A_{21} & A_{22} & \cdots & A_{2d} \\
\hline
\vdots & \vdots & \ddots & \vdots \\
\hline
A_{d1} & A_{d2} & \cdots & A_{dd}
\end{array}\right].
\]

\noindent
(1)$\Leftarrow $(2): Suppose the condition (1) holds. For each $i=1,\ldots,n,$ comparing the  $(i,i)$-blocks of $A\cdot {\rm diag}(P_{g_1},\ldots,P_{g_d})$ and  ${\rm diag}(P_{g_1},\ldots,P_{g_d})\cdot A$, we get
\[
A_{ii}P_{g_i}=P_{g_i}A_{ii}.
\]
This equality holds for every $g_i\in G_i$. Hence by Proposition \ref{prop:characterization_G_circulant}, $A_{ii}$ is a $G_i$-circulant matrix.

Now for each $1\leq i\not=j\leq n$, comparing the $(i,j)$-blocks of of $A\cdot {\rm diag}(P_{g_1},\ldots,P_{g_d})$ and  ${\rm diag}(P_{g_1},\ldots,P_{g_d})\cdot A$, we get
\[
A_{ij}P_{g_j}=P_{g_i}A_{ij}.
\]
This equality holds for every $g_i\in G_i$ and $g_j\in G_j$. In particular, $A_{ij}P_{g_j}=A_{ij}$, for every $g_j\in G_j$. This implies that all the columns of $A_{ij}$ are equal. Similarly, since $A_{ij}=P_{g_i}A_{ij}$, for every $g_i\in G_i$, all the rows of $A_{ij}$ are equal. Hence $A_{ij}$ of the form $\alpha_{ij }\ones$. We conclude that $A$ is in $\sJ_{G_1,\ldots,G_d}(R)$.

\medskip\noindent
(2)$\Rightarrow $(1): This is clear from a straightforward  multiplication of matrices and Proposition \ref{prop:characterization_G_circulant}.
\end{proof}

\begin{cor} 
If a matrix $A$ in $\sJ_{G_1,\ldots,G_d}(R)$ is invertible as an element of $M_{n,n}(R)$, then it is a unit in the ring $\sJ_{G_1,\ldots,G_d}(R)$.
\end{cor}
\begin{proof}
Since $A$ commutes with all ${\rm diag}(P_{g_1},\ldots,P_{g_d})$, so does its matrix inverse $A^{-1}$.
\end{proof}

\begin{prop} \label{prop:joingroupring}
 $\mathcal{J}_{G_1,\dots,G_d}(R)$ is a  subring of $M_n(R)$, and if $R$ is commutative, an $R$-subalgebra of $M_n(R)$.
\end{prop}

\begin{proof}
It is clear that if $A$ and $B$ in $ M_n(R)$ commute with ${\rm diag}(P_{g_1},\ldots,P_{g_d})$ then $A+B$ and $AB$ also commute with ${\rm diag}(P_{g_1},\ldots,P_{g_d})$.
\end{proof}
Here is another direct proof for Proposition \ref{prop:joingroupring}.
\begin{proof}
It is immediate that $\mathcal{J}_{G_1,\dots,G_d}(R)$ is closed with respect to linear combinations. The product of two matrices $A,B$ in $\sJ_{G_1,\dots,G_d}(R)$ can be performed blockwise:

\begin{equation}\label{eq:product}
A\cdot B=\left[\begin{array}{c|c|c|c}
C_1 & a_{12}\ones & \cdots & a_{1d}\ones \\
\hline
a_{21}\ones & C_2 & \cdots & a_{2d}\ones \\
\hline
\vdots & \vdots & \ddots & \vdots \\
\hline
a_{d1}\ones & a_{d2}\ones & \cdots & C_d
\end{array}\right]
\cdot
\left[\begin{array}{c|c|c|c}
D_1 & b_{12}\ones & \cdots & b_{1d}\ones \\
\hline
b_{21}\ones & D_2 & \cdots & b_{2d}\ones \\
\hline
\vdots & \vdots & \ddots & \vdots \\
\hline
b_{d1}\ones & b_{d2}\ones & \cdots & D_d
\end{array}\right].
\end{equation}

In the product matrix, a typical  diagonal entry $(AB)_{ll}$  is of the form
\[C_l D_l +\sum_{j=1,j\neq l}^dk_ja_{lj}b_{jl}\ones,\] and a typical off-diagonal $(AB)_{st}$ ($s\ne t $) is of the form 
\[a_{st}\ones D_t+\sum_{j \neq s, j\neq t}^{d}k_ja_{sj}b_{jt}\ones+C_s b_{st}\ones.
\]

 Since the product of $G$-circulant matrices is a $G$-circulant matrix (by Corollary \ref{cor:groupring}), the diagonal blocks of the product are $G$-circulant. Moreover, since $G$-circulant matrices are semimagic (Remark \ref{rem:semimagic}), the non-diagonal blocks of the product are constant matrices.  This shows that $\mathcal{J}_{G_1,\dots,G_d}(R)$ is closed with respect to products.
\end{proof}

\subsection{Augmentation Map} Recall that the augmentation map for a group ring $R[G]$ is the map $\epsilon \colon R[G] \rightarrow R$ defined  by $\epsilon(\sum \alpha_g g) = \sum_g \alpha_g$. By identifying $R[G]$ with $G$-circulant matrices, we can  
  generalize the augmentation map to joins of group rings. It is a map $\epsilon\colon \sJ_{G_1,\ldots,G_d}(R)\to M_d(R)$ defined as follows. For $A \in \sJ_{G_1,\ldots,G_d}(R)$ given by 
\[A=\begin{bmatrix} 
A_1&a_{12} J_{k_1,k_2} &\cdots & a_{1d}J_{k_1,k_d}\\
a_{21}J_{k_2,k_1} &A_2 &\cdots & a_{2d}J_{k_2,k_d}\\
\vdots&\vdots& &\vdots\\
a_{d1}J_{k_d,k_1}&a_{d2}J_{k_d,k_2}&\cdots& A_d
\end{bmatrix},\]
define $\epsilon (A)$ to be the $d \times d$ matrix obtained by replacing each block of $A$ with the corresponding row sum. That is, 
\[\epsilon (A) = \begin{bmatrix}
\epsilon(A_1) & k_2 a_{12} & \cdots  & k_d a_{1d} \\
k_1 a_{21} & \epsilon(A_2) & \cdots & k_d a_{2d} \\
\vdots&\vdots&  &\vdots\\
k_1a_{n1} & k_2a_{n2}& \cdots & \epsilon(A_d)
\end{bmatrix}.\]
As in the above display, we use $\epsilon$ to denote the augmentation map on the join ring and the group ring level; the correct interpretation should always be clear from the context.
 We now show that this augmentation map is a ring homomorphism. It turns out that this statement is true more generally for block matrices whose blocks have constant row sums. Note that $A_i$ being a $G$-circulant matrix, it has the property that the sum of the entries in any row  or column is $\epsilon(A_i)$ by Remark \ref{rem:semimagic}. To prove this more general statement, we begin with a lemma.

\begin{lem} \label{generalaugmentation}
Let $\mathcal{C}$ denote the collection of arbitrary matrices $X$ over a unital ring with the property that the sum of any row is constant, denoted $\epsilon(X)$. For any two matrices $A$ and $B$ in $\mathcal{C}$, we have: 
\begin{enumerate}
    \item If $A+B$ is defined, then  $\epsilon(A+B) = \epsilon(A)+ \epsilon(A)$\\
    \item If $AB$ is defined, then $\epsilon(AB) = \epsilon(A)\epsilon(B)$
\end{enumerate}
\end{lem}

\begin{proof}
(1) is obvious. For (2), we show that $AB$ has constant row sums. Let $A$ be an $m \times n$ matrix and $B$ be an $n \times q$, so the product $AB$ is defined. The sum of the entries in the $i$ th row of $AB$ is given by
\[
 \sum_{j = 1}^q (AB)_{ij} 
  =  \sum_{j=1}^q \sum_{k=1}^n A_{ik}B_{kj}
  =  \sum_{k=1}^n A_{ik}\left( \sum_{j=1}^q B_{kj} \right)
  =\sum_{k=1}^n A_{ik}\left( \epsilon(B) \right)\\
 = \epsilon(A) \epsilon(B).
\]
This shows that the row sums of $AB$ are constant, and  that $\epsilon(AB) =\epsilon(A)\epsilon(B)$.
\end{proof}

\begin{prop}
Let $A$ and $B$ be $d \times d$ block matrices where each block belongs to set $\mathcal{C}$. Then the map $\epsilon$ which sends each such block matrix to the $d \times d$ matrix obtained by replacing each block entry with the corresponding row sum respects addition and multiplication is a ring homomorphism. In particular, the augmentation map on the join ring is a ring homomorphism.
\end{prop}

\begin{proof}
The fact that $\epsilon$ respects addition is clear. It remains to show that $\epsilon(AB) = \epsilon(A) \epsilon(B)$.  We will show that the $(i, j)$th entry is the same on both sides.
\begin{eqnarray*}
(\epsilon(AB))_{ij} & = &  \epsilon( \sum_{k=1}^n A_{ik} B_{kj})\\
& = & \sum \epsilon(A_{ik}) \epsilon(B_{kj})  \ \ \text{(by  Lemma \ref{generalaugmentation})}\\ 
& =& \sum (\epsilon(A))_{ik} (\epsilon(B))_{kj} \\
& = & (\epsilon(A) \epsilon(B))_{ij}
\end{eqnarray*}

The last statement about the augmentation map now follows because elements in a join of group rings have the property that their blocks have constant row sums; see Remark \ref{rem:semimagic}.
\end{proof}

 \subsection{The center of $\sJ_{G_1,\ldots,G_d}(R)$} 
 
 \begin{lem}
 If a matrix $A$  is in  $\sJ_{G_1,\ldots,G_d}(R)$ then its transpose $A^T$ is also in $\sJ_{G_1,\ldots,G_d}(R)$.
 \end{lem}
 \begin{proof}
 This follows easily from Lemma~\ref{lem:transpose_circulant}.
 \end{proof}
 
For each $1\leq i\leq d$, let $E_{ii}$ be the matrix which has a 1 at the $(i,i)$-position and zeros in all other positions.

\begin{lem}
  If $A$ is in the center of $\sJ_{G_1,\ldots,G_d}(R)$ then both $\epsilon(A)$ and $\epsilon(A^T)$ commute with $E_{ii}$ for every $1\leq i\leq d$.
\end{lem}

\begin{proof} Let $B_i$ be the matrix in $\sJ_{G_1,\ldots,G_d}(R)$  defined as follows: all blocks  of $B_i$ are zeros except the $(i,i)$-diagonal block, which is the appropriate identity matrix. 
 Clearly $\epsilon(B_i)=\epsilon(B_i^T)=E_{ii}$, and $A$ commutes with $B_i$. Since $B_i$ is a $(0,1)$-matrix, one has $(XB_i)^T=B_i^TX^T$ and $(B_iY)^T=Y^TB_i^T$, for all matrices $X,Y$. Hence 
\[
A^T B_i^T=(B_iA)^T=(AB_i)^T=B_i^TA^T,
\] 
i.e., $A^T$ commutes with $B_i^T$. The lemma follows since $\epsilon$ is a ring homomorphism.
\end{proof}

\begin{prop} \label{prop:center}
The center of $\sJ_{G_1,\dots,G_d}(R)$ consists of the matrices \eqref{eq:join circulant matrix} such that all $a_{ij}=0$, all $C_i$ have the same row sum: $\epsilon(C_{1})=\dots=\epsilon(C_{d})$,  and each $C_i$ is in the center of $R[G_i]$. 
\end{prop}

\begin{proof}
Let 
\[A=\begin{bmatrix} 
C_1&a_{12} J_{k_1,k_2} &\cdots & a_{1d}J_{k_1,k_d}\\
a_{21}J_{k_2,k_1} &C_2 &\cdots & a_{2d}J_{k_2,k_d}\\
\vdots&\vdots& &\vdots\\
a_{d1}J_{k_d,k_1}&a_{d2}J_{k_d,k_2}&\cdots& C_d
\end{bmatrix}\] 
be in the center of $\sJ_{G_1,\dots,G_d}(R)$. Then by the previous lemma, 
\[\epsilon(A)
=\begin{bmatrix} 
\epsilon(C_1)&k_2a_{12} &\cdots & k_da_{1d}\\
k_1a_{21} &\epsilon(C_2) &\cdots & k_da_{2d}\\
\vdots&\vdots& &\vdots\\
k_1a_{d1}&k_2a_{d2}&\cdots& \epsilon(C_d)
\end{bmatrix}
\] commutes with $E_{ii}$, for every $i=1,\ldots,d$. 
This implies that $\epsilon(A)$ is diagonal. By applying the above argument to $A^T$, we see that 
\[\epsilon(A^T)
=\begin{bmatrix} 
\epsilon(C_1)&k_2a_{21} &\cdots & k_da_{d1}\\
k_1a_{12} &\epsilon(C_2) &\cdots & k_da_{d2}\\
\vdots&\vdots& &\vdots\\
k_1a_{1d}&k_2a_{2d}&\cdots& \epsilon(C_d)
\end{bmatrix}\]
 is also diagonal. Thus, $k_i a_{ij}=k_ja_{ij}=0$ for every $i\not=j$.

Now let 
\[B=\begin{bmatrix} 
D_1&b_{12} J_{k_1,k_2} &\cdots & b_{1d}J_{k_1,k_d}\\
b_{21}J_{k_2,k_1} &D_2 &\cdots & b_{2d}J_{k_2,k_d}\\
\vdots&\vdots& &\vdots\\
b_{d1}J_{k_d,k_1}&b_{d2}J_{k_d,k_2}&\cdots& D_d
\end{bmatrix}\]
be any matrix in $\sJ_{G_1,\dots,G_d}(R)$. The $(1,1)$-block of $AB$ is
\[
C_1D_1+k_2a_{12}b_{21}J_{k_1,k_1}+\cdots+k_da_{1d}b_{d1}J_{k_1,k_1}=C_1D_1.
\]
The $(1,1)$-block of $BA$ is 
\[
D_1C_1+b_{12}k_2a_{21}J_{k_1,k_1}+\cdots+b_{1d}k_da_{d1}J_{k_1,k_1}=D_1C_1.
\]
Hence $C_1D_1=D_1C_1$. This implies that $C_1$ is in the center of $R[G_1]$. Similarly $C_i$ is   in the center of $R[G_i]$.

Next, we compare the $(1,2)$-blocks of $AB$ and $BA$. The $(1,2)$-block of $AB$ is
\[
\begin{aligned}
\epsilon(C_1)b_{12}J_{k_1,k_2}+a_{12}\epsilon(D_2)J_{k_1,k_2}+k_3a_{13}b_{32}J_{k_1,k_2}+\cdots+k_da_{1d}b_{d2}J_{k_1,k_2}\\
=\epsilon(C_1)b_{12}J_{k_1,k_2}+a_{12}\epsilon(D_2)J_{k_1,k_2}.
\end{aligned}
\]
The $(1,2)$-block of $BA$ is
\[
\begin{aligned}
\epsilon(D_1)a_{12}J_{k_1,k_2}+b_{12}\epsilon(C_2)J_{k_1,k_2}+k_3b_{13}a_{32}J_{k_1,k_2}+\cdots+k_db_{1d}a_{d2}J_{k_1,k_2}\\
=\epsilon(D_1)a_{12}J_{k_1,k_2}+b_{12}\epsilon(C_2)J_{k_1,k_2}.
\end{aligned}
\]
Hence 
\[
\epsilon(C_1)b_{12}+a_{12}\epsilon(D_2)=\epsilon(D_1)a_{12}+b_{12}\epsilon(C_2).
\]
By choosing $B$ such that $b_{12}=0$, $\epsilon(D_1)=0$ and $\epsilon(D_2)=1$, we imply that $a_{12}=0$. Then by choosing $B$ such that $b_{12}=1$, we imply that $\epsilon(C_1)=\epsilon(C_2)$. 
Similarly, we obtain that $a_{ij}=0$ for every $i\not= j$ and $\epsilon(C_1)=\epsilon(C_2)=\cdots=\epsilon(C_d)$. 

 For the converse implication, suppose that $A$ is of the form
\[A=\begin{bmatrix} 
C_1&0 &\cdots & 0\\
0 &C_2 &\cdots & 0\\
\vdots&\vdots& &\vdots\\
0&0&\cdots& C_d
\end{bmatrix},\]
where each $C_i$ is in the center of $R[G_i]$ and $\epsilon(C_1)=\epsilon(C_2)=\cdots=\epsilon(C_d)$. Let
\[B=\begin{bmatrix} 
D_1&b_{12} J_{k_1,k_2} &\cdots & b_{1d}J_{k_1,k_d}\\
b_{21}J_{k2,k_1} &D_2 &\cdots & b_{2d}J_{k_2,k_d}\\
\vdots&\vdots& &\vdots\\
b_{d1}J_{k_d,k_1}&b_{d2}J_{k_d,k_2}&\cdots& D_d
\end{bmatrix}\]
be any matrix in $\sJ_{G_1,\dots,G_d}(R)$.
Then 
\[
\begin{aligned}
AB&=\begin{bmatrix} 
C_1D_1&C_1b_{12} J_{k_1,k_2} &\cdots & C_1b_{1d}J_{k_1,k_d}\\
C_2b_{21}J_{k_2,k_1} &C_2D_2 &\cdots & C_2b_{2d}J_{k_2,k_d}\\
\vdots&\vdots& &\vdots\\
C_db_{d1}J_{k_d,k_1}&C_db_{d2}J_{k_d,k_2}&\cdots& C_dD_d
\end{bmatrix}\\
&=\begin{bmatrix} 
C_1D_1&\epsilon(C_1)b_{12} J_{k_1,k_2} &\cdots & \epsilon(C_1)b_{1d}J_{k_1,k_d}\\
\epsilon(C_2)b_{21}J_{k_2,k_1} &C_2D_2 &\cdots & \epsilon(C_2)b_{2d}J_{k_2,k_d}\\
\vdots&\vdots& &\vdots\\
\epsilon(C_d)b_{d1}J_{k_d,k_1}&\epsilon(C_d)b_{d2}J_{k_d,k_2}&\cdots& C_dD_d
\end{bmatrix},
\end{aligned}
\]
and using the fact that the matrices $C_i$ are  semimagic, we get
\[
\begin{aligned}
BA&=\begin{bmatrix} 
D_1C_1&b_{12} J_{k_1,k_2}C_2 &\cdots & b_{1d}J_{k_1,k_d}C_d\\
b_{21}J_{k_2,k_1}C_1 &C_2D_2 &\cdots & b_{2d}J_{k_2,k_d}C_d\\
\vdots&\vdots& &\vdots\\
b_{d1}J_{k_d,k_1}C_1&b_{d2}J_{k_d,k_2}C_2&\cdots& D_dC_d
\end{bmatrix}\\
&=\begin{bmatrix} 
D_1C_1&b_{12}\epsilon(C_2) J_{k_1,k_2} &\cdots & b_{1d}\epsilon(C_d)J_{k_1,k_d}\\
b_{21}J_{k_2,k_1}\epsilon(C_1) &C_2D_2 &\cdots & b_{2d}J_{k_2,k_d}\epsilon(C_d)\\
\vdots&\vdots& &\vdots\\
b_{d1}J_{k_d,k_1}\epsilon(C_1)&b_{d2}J_{k_d,k_2}\epsilon(C_2)&\cdots& D_dC_d
\end{bmatrix}
\end{aligned}
\]
Clearly (using the fact that  $\epsilon$ maps the center of $R[G_i]$ into the center of $R$), $AB=BA$. Hence $A$ is in the center of $\sJ_{G_1,\dots,G_d}(R)$.

\end{proof}
\begin{cor} \label{cor:dimension_center}
For a join algebra $\sJ=\mathcal{J}_{G_1,\dots,G_d}(k)$ where $k$ is a field, the $k$-dimension of the center of $\sJ$ is  $\sum_{i=1}^d \dim_k(Z(k[G_i]))-d+1=\sum_{i=1}^dc(G_i)-d+1$ where $Z(k[G_i])$ denotes the center of the group algebra and $c(G_i)$ denotes the number of conjugacy classes in $G_i$. \end{cor}
We will shortly see (Proposition \ref{prop:number_of_irr_modules}) that in the case where $\sJ$ is semisimple and $k$ is algebraically closed, $\sum_{i=1}^dc(G_i)-d+1$ also counts the number of irreducible modules over $\sJ$.

\subsection{Some further ring-theoretic properties of $\mathcal{J}_{G_1, G_2, \ldots, G_d}(R)$}
Next, we discuss some further ring-theoretic properties of $\mathcal{J}_{G_1, G_2, \ldots, G_d}(R)$. Before doing so, we first recall some notations in ring theory. Our main reference for this part is \cite{lam}. We first recall the definition of semisimplicity.

\begin{defn}(\cite[Theorem and Definition 2.5]{lam})

Let $M$ be a module over a ring $R$. Then 
\begin{enumerate}
\item We say that $M$ is semisimple if every $R$-submodule of $M$ is an $R$-submodule direct summand of $M.$
 \item We say that $R$ is semisimple if all left $R$-modules are semisimple. Equivalently, $R$ is semisimple if the left regular $R$-module $_{R}R$ is semisimple. 
 \end{enumerate} 
\end{defn}

A closely related notion of semisimplicity is Jacobson semisimplicity (or $J$-semisimplicity), which we now recall. 
\begin{defn}
Let $R$ be a ring. The Jacobson radical of $R$, denoted by $\Rad(R)$, is the intersection of all maximal left ideals in $R.$ The ring $R$ is called Jacobson semisimple if $\Rad(R)=0.$
\end{defn}
A famous theorem of Maschke says the following.
\begin{thm}(Maschke's theorem, \cite[Theorem 6.1]{lam}) \label{thm:classical_maschke}

Let $R$ be a semisimple ring and $G$ a finite group such that $|G|$ is invertible in $R$. Then the group algebra $R[G]$ is semisimple. \end{thm}

The following theorem is a natural generalization of Maschke's theorem to  $\mathcal{J}_{G_1, G_2, \ldots, G_d}(R)$. 

\begin{thm} (Generalized Maschke's theorem) \label{thm:maschke}
Let $R$ be a semisimple ring. Suppose that $|G_i|$ is invertible in $R$ for all $1 \leq i \leq d.$ Then the ring $\sJ_{G_1, G_2, \ldots, G_d}(R)$ is semisimple. 
\end{thm}

In this section, we provide the first proof for this theorem. A 
second proof in the case that $R=k$ is a field can be found in Section $5$ where we explicitly describe the Jacobson radical of $J_{G_1, G_2, \ldots, G_d}(k)$.
To explain the first proof, we first discuss the structure of the group ring $R[G]$ when $|G|$ is invertible in $R.$ Let $\epsilon$ be the augmentation map $\epsilon\colon R[G] \to R.$ Let $\Delta(G) =\Delta_{R}(G)$ be the kernel of $\epsilon.$ Since $|G|$ is invertible in $R$, we can consider 
\[ e_{G}:=\frac{1}{|G|} \sum_{g \in G} g .\] 
We can easily see that $e_G$ is a central idempotent in $R[G]$. Furthermore, by \cite[Proposition 3.6.7]{GR}, we have 
\begin{prop} \label{prop:splitting}
Suppose that $|G|$ is invertible in $R$. Then $R[G]$ is a direct product of the following rings 
\[ R[G] \cong R[G]e_G \times R[G](1-e_G) .\] 
Furthermore 
\[ R[G] e_G \cong R, \]
and 
\[ R[G](1-e_G) = \Delta_{R}(G) .\] 
In particular, $\Delta(G)$ is a semisimple ring.
\end{prop}

A direct generalization of this structure theorem is the following. 
\begin{thm} \label{thm:decomposition}
Assume that $|G_i|$ is invertible in $R$ for all $1 \leq i \leq d.$ Then  $\sJ_{G_1, G_2, \ldots, G_d}(R)$ is a direct product of the following rings 
\[ \sJ_{G_1, G_2, \ldots, G_d}(R) \cong M_d(R) \times \prod_{i=1}^d \Delta_{R}(G_i) .\] 
\end{thm}

\begin{proof}
Let $f_i=f_{G_i}=1-e_{G_i} \in R[G_i]$ where $e_{G_i}$ is defined as above. Since the ring of all circulant matrices is isomorphic to the group ring $R[G_i]$, we can also consider $f_i$ as a $G_i$-circulant matrix. Let $\tilde{f}_i$ be the following matrix in $\sJ_{G_1, G_2, \ldots, G_d}(R)$
\[
\tilde{f}_i =\left[\begin{array}{c|c|c|c}
0 & 0 & \dots & 0 \\
\hline
0 & 0 & \dots & 0 \\
\hline
\vdots & \vdots & f_i & \vdots\\
\hline
0 & 0 & \dots & 0
\end{array}
\right].
\]
In other words, all blocks of $\tilde{f}_i$, except the $i$-diagonal block which is $f_i$,  are $0.$  Additionally, we define 
\[ \tilde{f}_{d+1}=I_n-\sum_{i=1}^d \tilde{f}_i=\bigoplus_{i=1}^{d} e_{G_i}. \] 
By definition, we have 
\[ \tilde{f}_i^2 = \tilde{f}_i, \forall 1 \leq i \leq d+1.\]
Furthermore, by Proposition \ref{prop:center}, $\tilde{f_i}$ belongs to the center of $\sJ_{G_1, G_2, \ldots, G_d}(R).$ We conclude that for all $i$, $\tilde{f}_i$ is a central idempotent in $\sJ_{G_1, G_2, \ldots,G_d}(R).$ Consequently, we have the following ring isomorphism 
\[ \sJ_{G_1, \ldots, G_d}(R) \cong \tilde{f}_{d+1} J_{G_1, \ldots, G_d}(R) \times \prod_{i=1}^d \tilde{f}_{i} J_{G_1, \ldots, G_d}(R) .\] 
Direct calculations show that for $1 \leq i \leq d$ 
\[ \tilde{f}_{i} \sJ_{G_1, \ldots, G_d}(R) \cong \Delta_R(G_i).\] We claim that the augmentation map 
\[ \epsilon\colon \sJ_{G_1, G_2, \ldots, G_d}(R) \to M_d(R) \] 
induces a ring isomorphism 
\[ \epsilon\colon \tilde{f}_{d+1} \sJ_{G_1, G_2, \ldots, G_d}(R) \to M_d(R).\] 
Since $\epsilon(\tilde{f}_{d+1})=1$, we have for every $A \in  \sJ_{G_1, G_2, \ldots, G_d}(R)$
\[ \epsilon(\tilde{f}_{d+1} A) = \epsilon(A) .\] 
Because the map $\epsilon\colon \sJ_{G_1, G_2, \ldots, G_d}(R) \to M_d(R)$ is surjective, the above equality shows that the induced map on $\tilde{f}_{d+1} \sJ_{G_1, G_2, \ldots, G_d}(R)$ is surjective as well. We claim that it is injective as well. In fact, suppose that $\tilde{f}_{d+1} A \in \ker(\epsilon)$ for some $A \in \sJ_{G_1, G_2, \ldots, G_d}(R).$ Suppose that $A$ has the following form 
\[A=\begin{bmatrix} 
C_1&a_{12} J_{k_1,k_2} &\cdots & a_{1d}J_{k_1,k_d}\\
a_{21}J_{k_2,k_1} &C_2 &\cdots & a_{2d}J_{k_2,k_d}\\
\vdots&\vdots& &\vdots\\
a_{d1}J_{k_d,k_1}&a_{d2}J_{k_d,k_2}&\cdots& C_d
\end{bmatrix}.\]
Then 
\[ 0= \epsilon(A) 
=\begin{bmatrix} 
\epsilon(C_1)&k_2a_{12} &\cdots & k_da_{1d}\\
k_1a_{21} &\epsilon(C_2) &\cdots & k_da_{2d}\\
\vdots&\vdots& &\vdots\\
k_1a_{d1}&k_2a_{d2}&\cdots& \epsilon(C_d)
\end{bmatrix}.
\]
Since $k_i=|G_i|$ are invertible in $R$, we conclude that $a_{ij}=0.$ Additionally 
\[ \epsilon(C_1) = \epsilon(C_2) = \cdots =\epsilon(C_d) = 0 .\] 
We then see that $\tilde{f}_{d+1}A=0.$ We conclude that the   map 
\[ \epsilon: \tilde{f}_{d+1} \sJ_{G_1, G_2, \ldots, G_d}(R) \to M_d(R),\] 
is injective. 
\end{proof}

We are ready to prove Theorem \ref{thm:maschke}. 

\medskip
\begin{proof}[\textbf{Proof of Theorem \ref{thm:maschke}}] Since $R$ is semisimple, by Morita equivalence, we conclude that $M_d(R)$ is semisimple. Furthermore, by the classical Maschke theorem \ref{thm:classical_maschke} and the decomposition mentioned in Proposition \ref{prop:splitting}, we know that $\Delta(G_i)$ is semisimple for $1 \leq i \leq d.$  As a product of semisimple rings is semisimple, Theorem \ref{thm:decomposition} implies that $\sJ_{G_1, G_2, \ldots, G_d}(R)$ is semisimple as well. 
\end{proof}

We discuss some consequences of Theorem \ref{thm:decomposition}. 

\begin{prop} \label{prop:number_of_irr_modules}
Suppose that $k$ is algebraically closed and ${\rm char}(k)$ is relatively prime to $\prod_{i=1}^d |G_i|$. Then the number of irreducible modules over $\sJ_{G_1, G_2, \ldots, G_d}(k)$ is 
\[ c(G_1)+c(G_2)+\cdots+c(G_d)-d+1 .\] 
where $c(G_i)$ is the number of conjugacy classes of $G_i$.
\end{prop}

\begin{proof}
By Artin-Wedderburn theorem, the number of simple modules over $k[G_i]$ is $c(G_i).$ From the decomposition $k[G_i] \cong k \times \Delta_{k}(G)$, we conclude that the number of irreducible modules over $\Delta_{k}(G)$ is $c(G_i)-1.$ Finally, by the Morita equivalence, there is exactly one irreducible module over $M_d(k)$. Therefore, by Theorem \ref{thm:decomposition}, we conclude that the number of irreducible modules over $\sJ_{G_1, G_2, \ldots, G_d}(k)$ is 
\[ 1+ \sum_{i=1}^d (c(G_i)-1) =\sum_{i=1}^d c(G_i)-d+1.
\qedhere\] 
\end{proof}

Let us discuss some special cases.

\begin{ex} \label{ex:cyclic}
Suppose that $G_i=\Z/k_i$, the cyclic group of order $k_i.$ Let $R=k$ be an algebraically closed field of characteristic $0$. Then, we have 
\[ k[G_i] \cong k[x]/(x^{k_i}-1) \cong \prod_{j=1}^{k_i} k .\]
By Theorem \ref{thm:decomposition}, we conclude that 
\[ \sJ_{G_1, G_2, \ldots, G_d}(k) \cong M_d(k) \times k^{n-d},\]
where $n =\sum_{i=1}^d k_i.$ Please see the last section for an explicit map for this isomorphism. 
\end{ex}

\begin{ex} \label{ex:join}
Let us consider $G_1=D_{2n}$,  $G_2= \Z/2$. Suppose that $k$ is an algebraically closed field of characteristic $0.$ By the  Artin-Wedderburn theorem (see \cite[Section 18.3]{james2001representations}), we know that 
\[ k[D_{2n}] \cong  \begin{cases}
  k^2 \times M_2(k)^{\frac{n-1}{2}}  & \text{if } n \text{ is odd} \\
  k^4 \times M_2(k)^{\frac{n-2}{2}} & \text{if } n \text{ is even}.
\end{cases} \] 
Consequently 
\[ \Delta(D_{2n})\cong  \begin{cases}
  k \times M_2(k)^{\frac{n-1}{2}}  & \text{if } n \text{ is odd} \\
  k^3 \times M_2(k)^{\frac{n-2}{2}} & \text{if } n \text{ is even}.
\end{cases} \] 
We also know that 
\[ \Delta(\Z/2)= k.\] 
Consequently, by Theorem \ref{thm:decomposition}, we have 
\[ \sJ_{D_{2n}, \Z/2}(k)  \cong  \begin{cases}
  k^2 \times M_2(k)^{\frac{n+1}{2}}  & \text{if } n \text{ is odd} \\
  k^4 \times M_2(k)^{\frac{n}{2}} & \text{if } n \text{ is even}. 
\end{cases} \] 
We conclude that 
\[ \sJ_{D_{2n}, \Z/2}(k) \cong k[D_{2(n+2)}] .\] 
\end{ex}
We give another example of $G$ such that the group algebra of  $k[G]$ is isomorphic to the join algebra $\sJ_{G_1, G_2, \ldots, G_d}(k)$ for some $d \geq 2.$ 
\begin{ex}  \label{ex:join1}
Assume that $k$ is an algebraically closed field whose characteristic differernt from $2$, $3$ and $7$. Then the group algebra $k[S_4]$ is isomorphic to any join algebra $\sJ_{G_1, G_2}(k)$ where $G_1$ is the trivial group and $G_2$ is the group of order $21$ with the following presentation 
\[ \langle x,y: x^7=y^3=1, yxy^{-1}=x^2 \rangle .\] 
In fact, from the representations of $G_2$ (see \cite[Theorem 25.10]{james2001representations}) we have 
\[ k[G_2] \cong k^{3} \times M_3(k)^2.  \] 
Additionally, from the representations of $S_4$ (see \cite[Example 16.3]{james2001representations}), we have 

\[ k[S_4] \cong k^2 \times M_2(k) \times M_3(k)^2 .\] 
By Theorem \ref{thm:decomposition} we have 
\[ \sJ_{G_1, G_2}(k) \cong k^2 \times M_2(k) \times M_3(k)^2 .\]
Consequently 
\[ k[S_4] \cong \sJ_{G_1, G_2}(k) .\]
\end{ex}
\begin{rem}
Let $k$ be a field. By Corollary \ref{cor:dimension_center}, the join algebra $\sJ_{G_1, G_2, \ldots, G_d}(k)$ is not abelian when $d \geq 2.$ Therefore, it cannot be isomorphic to $k[G]$ where $G$ is an abelian group. In Example \ref{ex:join} and Example \ref{ex:join1}, we have provided two examples of non-abelian groups such that their group rings are isomorphic to a join ring $\sJ_{G_1, G_2, \ldots, G_d}(k)$ with $d \geq 2.$ We can show also that the group algebra of $A_5$, the alternating group on $5$ letters, is not isomorphic to $\sJ_{G_1, G_2, \ldots, G_d}(k)$ for $d \geq 2.$ We wonder whether the same statement is true for other non-abelian simple groups?
\end{rem}

\section{Characterization of units in the algebra $\mathcal{J}_{G_1, G_2, \ldots, G_d}(k)$}

\bigskip

Let finite groups $G_{1},G_{2},\ldots ,G_{d}$ be of
respective orders $k_{1},k_{2},\ldots ,k_{d}$. We characterize the units in
the join $\mathcal{J}_{G_1, G_2, \ldots, G_d}(k)$, where $k$ is any field.

Note that $\mathcal{J}_{G_1, G_2, \ldots, G_d}(k)$ is a subalgebra of $M_n(k)$, where $n = k_1+k_2 +  \cdots + k_d$. We have shown that an element $X$ in $\mathcal{J}_{G_1, G_2, \ldots, G_d}(k)$ is a unit in $\mathcal{J}_{G_1, G_2, \ldots, G_d}(k)$  if and only if it is an invertible matrix in $M_n(k)$. Therefore, characterizing units in $\mathcal{J}_{G_1, G_2, \ldots, G_d}(k)$ reduces to determining necessary and sufficient conditions under which an element in the join is an invertible matrix. To this end, we need the following definition.

\begin{defn}  A $G_i$-circulant matrix $A$ is said to be almost invertible if and only if $\epsilon(A)=0$ and nullity$(A)=1$.   \end{defn}   \ The result is:

\bigskip

\begin{thm} \label{unitstheorem}
An element \[
X=\left[\begin{array}{c|c|c|c}
C_1 & a_{12}\ones & \cdots & a_{1d}\ones \\
\hline
a_{21}\ones & C_2 & \cdots & a_{2d}\ones \\
\hline
\vdots & \vdots & \ddots & \vdots \\
\hline
a_{d1}\ones & a_{d2}\ones & \cdots & C_d
\end{array}\right],
\] is a unit in $\mathcal{J}_{G_1, G_2, \ldots, G_d}(k)$ if and only if the ``principal diagonal matrices" $C_1,C_2,\ldots,C_d$ are each (independently) either
invertible or almost invertible and $\epsilon(X)$ is an invertible matrix. Here $\epsilon$ is the augmentation map $\sJ_{G_1, G_2, \ldots, G_d}(k) \to M_d(k)$.
\end{thm}

\begin{proof}  To keep notation simple, we give the proof here for $d=3$\, writing $m=k_{1},n=k_{2},q=$ $k_{3}$ and setting $X:=\left[ 
\begin{array}{lll}
A & \alpha \ones_{m,n} & \beta \ones_{m,q} \\ 
\gamma \ones_{n,m} & B & \delta \ones_{n,q} \\ 
\nu \ones_{q,m} & \eta \ones_{q,n} & C%
\end{array}%
\right] $,  but
the argument goes through for arbitrary $d$. (The $\ones$  matrices are subscripted with their dimensions for convenient reference.) First, let us observe that the conditions that $A, B, C$ are each
(independently) either invertible or almost invertible, and that $\epsilon(X)$ is
an invertible matrix, are necessary. \ The latter condition is clearly
necessary as $\epsilon$ here is a ring homomorphism.  Now take $A$ for example
and suppose $A$ is not invertible, but that $X$ is a unit.  Then if $A$ has
nullity $2$ or more, we can do row operations on $A$ to produce two rows of
zeros; those operations on the first $m$ rows of $X$ will produce a matrix
of the form $\left[ 
\begin{array}{lll}
A_{1} & \alpha \ones_{1,n} & \beta \ones_{1,q} \\ 
\vdots & \vdots & \vdots \\ 
A_{m-2} & \alpha \ones_{1,n} & \beta \ones_{1,q} \\ 
\overrightarrow{0} & r\alpha \ones_{1,n} & r\beta \ones_{1,q} \\ 
\overrightarrow{0} & s\alpha \ones_{1,n} & s\beta \ones_{1,q}%
\end{array}%
\right] $\ (where we may have permuted rows of $A$) in which the last two
rows \ are linearly dependent, a contradiction; so the nullity of $A$ is $1$.  Further, if $r_{1}A_{1}+\cdots +r_{m}A_{m}=0$ is a relation on the rows
of $A$ with, say, $r_{m}\neq 0$, then summing all entries gives $\epsilon(A)\sum
r_{i}=0.$ Now if $\epsilon(A)\neq 0$ (so $\sum r_{i}=0$) we can do elementary row
operations to replace $A_{m}$ with $r_{1}A_{1}+\cdots +r_{m}A_{m}$; the same
operations performed on the first $m$ rows of $X$ produce a row of $0$'s. 
So we must have $\epsilon(A)=0$, and similarly for $B, C$.

\bigskip

Now suppose the necessary conditions are in place.    Suppose a linear
combination of the rows of $X$ with coefficients $r_{1},\cdots
,r_{m},s_{1},\cdots ,s_{n},t_{1},\cdots ,t_{q}$ is the zero vector.  Let $
r=\sum r_{i},s=\sum s_{i},t=\sum t_{i}$.  Then, considering the first $m$
columns, then the next $n$ columns, and finally the last $q$ columns of $X$
we have the equations

\begin{eqnarray*}
\sum r_{i}A_{i}+s\gamma J_{1,m}+t\nu J_{1,m} &=&\overrightarrow{0}%
_{m} \\
r\alpha J_{1,n}+\sum s_{i}B_{i}+t\eta J_{1,n} &=&\overrightarrow{0}_{n} \\
r\beta J_{1,q}+s\delta J_{1,q}+\sum t_{i}C_{i} &=&\overrightarrow{0}_{q}.
\end{eqnarray*}

Summing all entries in these vector equalities, we get%
\begin{eqnarray*}
r\epsilon(A)+s\gamma m+t\nu m &=&0 \\
r\alpha n+s\epsilon(B)+t\eta n &=&0 \\
r\beta q+s\delta q+t\epsilon(C) &=&0
\end{eqnarray*}
or 
\[
\left( 
\begin{array}{lll}
r & s & t%
\end{array}%
\right) \epsilon(X)=\left( 
\begin{array}{lll}
0 & 0 & 0%
\end{array}%
\right) 
\]%
so by invertibility of $\epsilon(X)$ we have $r=s=t=0$. \ Returning to our
displayed set of vector equations, we now have

\bigskip 
\begin{eqnarray*}
\sum r_{i}A_{i} &=&\overrightarrow{0}_{m} \\
\sum s_{i}B_{i} &=&\overrightarrow{0}_{n} \\
\sum t_{i}C_{i} &=&\overrightarrow{0}_{q}.
\end{eqnarray*}

If $A$ is invertible, the first of these conditions forces all $r_{i}=0
$. \ But this also must be the case if $A$ is almost invertible; this
follows since any linear relation on the rows of $A$ must have $%
r_{1}=r_{2}=\cdots =r_{m}$ (recall that $\epsilon(A) = 0$ implies that the row of all 1's is an eigenvector for $A$ with eigenvalue $0$), from which $mr_{1}=\sum r_{i}=r=0$.  But having already $\epsilon(A)=0$ and $\epsilon(X)$ invertible, 
we cannot have $m=0\in k$ (lest we have a column of zeros in $\epsilon(X)$), and so $0=r_{1}=r_{2}=\cdots =r_{m}$. \
Similar arguments for $B$ and $C$ imply that all coefficients $r_{1},\ldots
,r_{m},s_{1},\ldots ,s_{n},t_{1},\ldots ,t_{q}$ must be $0$, so $X$ is in
fact invertible. \end{proof}

\subsection{Structure of the unit group of $\sJ_{G_1, G_2, \ldots, G_d}(k)$ in the modular case.} 
In this subsection, we will let $R$ denote a commutative ring such that $|G_i|=0$ in $R$ for all $i$.  When $R =k$ is field, this means that $|G_i|$ is divisible by the characteristic of $k$ for all $i$. This is the modular case, in which we can give a better characterization of units in $\sJ_{G_1, G_2, \ldots, G_d}(k)$.

For simplicity, we will denote  $\sJ_{G_1, G_2, \ldots, G_d}(R)$ by $\sJ$. Recall that we have the augmentation map, which is a ring homomorphism 
\[ \epsilon: \sJ \to M_d(R) ,\]
sending 
\[A=\begin{bmatrix} 
C_1&a_{12} J_{k_1,k_2} &\cdots & a_{1d}J_{k_1,k_d}\\
a_{21}J_{k_2,k_1} &C_2 &\cdots & a_{2d}J_{k_2,k_d}\\
\vdots&\vdots& &\vdots\\
a_{d1}J_{k_d,k_1}&a_{d2}J_{k_d,k_2}&\cdots& C_d
\end{bmatrix} \mapsto \begin{bmatrix} 
\epsilon(C_1)&k_2a_{12} &\cdots & k_da_{1d}\\
k_1a_{21} &\epsilon(C_2) &\cdots & k_da_{2d}\\
\vdots&\vdots& &\vdots\\
k_1a_{d1}&k_2a_{d2}&\cdots& \epsilon(C_d)
\end{bmatrix}.
\]
In fact, under the assumption that $k_i=|G_i|=0$ in $R$, we conclude that the map $\epsilon$ lands in the subset of diagonal matrices in $M_d(R).$ Consequently, $\epsilon$ induces a group homomorphism 
\[ \epsilon: \sJ^{\times} \to (R^{\times})^d ,\]
sending 
\[ A \mapsto (\epsilon(C_1), \epsilon(C_2), \ldots, \epsilon(C_d)) .\] 
We can see that this map is surjective. Let $U_1(\sJ)$ be its kernel (we can think of an element of $U_1(\sJ)$ as an analog of a principal unit in the classical group ring $R[G]$). Then, we have the following short exact sequence. 
\begin{equation} \label{eq:unit}
1 \to U_1(\sJ) \to \sJ^{\times} \to (R^{\times})^d \to 1. 
\end{equation}
 We can observe that there is a natural section $(R^{\times})^d \to \sJ^{\times}$ sending 
\[ (a_1, a_2, \ldots, a_d) \mapsto \begin{bmatrix} 
a_1 I_{k_1}& 0 &\cdots & 0\\
0 &a_2I_{k_2} &\cdots & 0\\
\vdots&\vdots& &\vdots\\
0&0 &\cdots& a_d I_{k_d}
\end{bmatrix}. \] 

Consequently, we have the following proposition.

\begin{prop} \label{prop:unit_decomposition}
$\sJ^{\times}$ is a semidirect product of $U_1(\sJ)$ and $(R^{\times})^d$:
\[ \sJ^{\times} \cong U_1(\sJ) \rtimes (R^{\times})^d .\] 
\end{prop}
Next, we investigate further the structure of $U_1(\sJ).$ Let $M=R^{d^2-d}$ be the abelian group of all $d \times d$ matrices  of the form 

\[ \begin{bmatrix}
0 & a_{12} & \cdots & a_{1d} \\
a_{21} & 0 & \cdots & a_{2d} \\
\vdots  & \vdots  & \ddots & \vdots  \\
a_{d1} & a_{d2} & \cdots & 0 \end{bmatrix},\] 
where the group structure is given by the usual component-wise matrix addition. 

We have the following observation.
\begin{prop}
The logarithm map 
\[ \log: U_1(\sJ) \to M ,\]
sending
\[
A=\left[\begin{array}{c|c|c|c}
A_1 & a_{12}\ones & \cdots & a_{1d}\ones \\
\hline
a_{21}\ones & A_2 & \cdots & a_{2d}\ones \\
\hline
\vdots & \vdots & \ddots & \vdots \\
\hline
a_{d1}\ones & a_{d2}\ones & \cdots & A_d
\end{array}\right]
\mapsto 
 M= \begin{bmatrix}
0 & a_{12} & \cdots & a_{1n} \\
a_{21} & 0 & \cdots & a_{2n} \\
\vdots  & \vdots  & \ddots & \vdots  \\
a_{n1} & a_{n2} & \cdots & 0 \end{bmatrix} ,\]
is a surjective group homomorphism. Furthermore, $\log$ 
has a left inverse $\psi: M \to U_1(\sJ)$ that sends $ M$ to 
\begin{equation*}
\psi(M)=\left[\begin{array}{c|c|c|c}
I & a_{12}\ones & \cdots & a_{1d}\ones \\
\hline
a_{21}\ones & I & \cdots & a_{2d}\ones \\
\hline
\vdots & \vdots & \ddots & \vdots \\
\hline
a_{d1}\ones & a_{d2}\ones & \cdots & I
\end{array}\right],
\end{equation*}
\end{prop}
We call the map in the statement of this proposition "log" because the matrix operations on the domain and codomain of this map are multiplication and addition, respectively.
\begin{proof}
Suppose 
\[ A=\left[\begin{array}{c|c|c|c}
A_1 & a_{12}\ones & \cdots & a_{1d}\ones \\
\hline
a_{21}\ones & A_2 & \cdots & a_{2d}\ones \\
\hline
\vdots & \vdots & \ddots & \vdots \\
\hline
a_{d1}\ones & a_{d2}\ones & \cdots & A_d
\end{array}\right]
\text{ and }  
B=\left[\begin{array}{c|c|c|c}
B_1 & b_{12}\ones & \cdots & b_{1d}\ones \\
\hline
b_{21}\ones & B_2 & \cdots & b_{2d}\ones \\
\hline
\vdots & \vdots & \ddots & \vdots \\
\hline
b_{d1}\ones & a_{d2}\ones & \cdots & B_d
\end{array}\right].
\]
The condition $\epsilon(A_i)=\epsilon(B_i)=1$ implies that 
\begin{equation} \label{eq:multiplicationAB}
AB=\left[\begin{array}{c|c|c|c}
A_1B_1 & (a_{12}+b_{12})\ones & \cdots & (a_{1d}+b_{1d})\ones \\
\hline
(a_{21}+b_{21})\ones & A_2B_2 & \cdots & (a_{2d}+b_{2d}) \ones \\
\hline
\vdots & \vdots & \ddots & \vdots \\
\hline
(a_{d1}+b_{d1}) \ones & (a_{d2}+b_{d2})\ones & \cdots & A_d B_d
\end{array}\right],
\end{equation}
This calculation shows that $\log$ is a group homomorphism. A similar calculation shows that $\psi$ is a group homomorphism and $\log \circ\, \psi$ is the identity map on $M.$
\end{proof}
Let $DU_1(\sJ) = \ker(\log)$. By definition, we can see that 
\[ DU_1(\sJ) \cong U_1(R[G_1]) \times U_1(R[G_2]) \times \ldots \times U_1(R[G_d]) ,\]
where $U_1(R[G_i])$ is the group of principal units in $R[G_i].$ We also have a short exact sequence 
\begin{equation} \label{eq:unit2}
1 \to DU_1(\sJ) \xrightarrow{\iota} U_1(\sJ) \to M \to 1.     
\end{equation}
It turns out that this exact sequence splits as well. 
\begin{prop}  \label{prop:principal_unit_decomposition} 
The short exact sequence \ref{eq:unit2} splits. In other words, 
\[ U_1(\sJ) \cong DU_1(\sJ) \times M \cong U_1(R[G_1]) \times U_1(R[G_2]) \times \ldots \times U_1(R[G_d]) \times R^{d^2-d}.\]
\end{prop}
\begin{proof}
Let us construct an inverse $\Phi$ of $\iota.$ Let $A$ be an element in $U_1(\sJ)$. Suppose that 

\begin{equation*}
A=\left[\begin{array}{c|c|c|c}
A_1 & a_{12}\ones & \cdots & a_{1d}\ones \\
\hline
a_{21}\ones & A_2 & \cdots & a_{2d}\ones \\
\hline
\vdots & \vdots & \ddots & \vdots \\
\hline
a_{d1}\ones & a_{d2}\ones & \cdots & A_d
\end{array}\right].
\end{equation*}
We define 
\begin{equation*}
\Phi(A)=\left[\begin{array}{c|c|c|c}
A_1 & 0 & \cdots & 0 \\
\hline
0 & A_2 & \cdots & 0 \\
\hline
\vdots & \vdots & \ddots & \vdots \\
\hline
0 & 0 & \cdots & A_d
\end{array}\right].
\end{equation*}
It is clear that $\Phi(\iota(A))=A$ for all $A \in DU_1(\sJ).$ We claim that $\Phi\colon  U_1(\sJ) \to DU_1(\sJ)$ is a group homomorphism. In fact, let $B$ be another element in $U_1(\sJ)$ 

\begin{equation*}
B=\left[\begin{array}{c|c|c|c}
B_1 & b_{12}\ones & \cdots & b_{1d}\ones \\
\hline
b_{21}\ones & B_2 & \cdots & b_{2d}\ones \\
\hline
\vdots & \vdots & \ddots & \vdots \\
\hline
b_{d1}\ones & a_{d2}\ones & \cdots & B_d
\end{array}\right].
\end{equation*}
We have 
\begin{equation*}
AB=\left[\begin{array}{c|c|c|c}
A_1B_1 & (a_{12}+b_{12})\ones & \cdots & (a_{1d}+b_{1d})\ones \\
\hline
(a_{21}+b_{21})\ones & A_2B_2 & \cdots & (a_{2d}+b_{2d}) \ones \\
\hline
\vdots & \vdots & \ddots & \vdots \\
\hline
(a_{d1}+b_{d1}) \ones & (a_{d2}+b_{d2})\ones & \cdots & A_d B_d
\end{array}\right].
\end{equation*}
From this equation, we can see that 
\[ \Phi(AB)= \Phi(A) \Phi(B) . \]
This shows that $\Phi$ is a group homomorphism, as required. 
\end{proof}
In summary, we have the following theorem about the structure of $\sJ^{\times}.$
\begin{thm} \label{units-modular}
\[ \sJ^{\times} \cong (R^{d^2-d} \times \prod_{i=1}^d U_1(R[G_i])) \rtimes (R^{\times})^d .\]
\end{thm}

We remark that the proof of Proposition \ref{prop:principal_unit_decomposition} shows a little more. 
\begin{cor}
Let $A$ be an element in $\sJ_{G_1, G_2, \ldots, G_d}(R)$ 
\begin{equation*}
A=\left[\begin{array}{c|c|c|c}
A_1 & a_{12}\ones & \cdots & a_{1d}\ones \\
\hline
a_{21}\ones & A_2 & \cdots & a_{2d}\ones \\
\hline
\vdots & \vdots & \ddots & \vdots \\
\hline
a_{d1}\ones & a_{d2}\ones & \cdots & A_d
\end{array}\right].
\end{equation*}
Then $A$ is invertible if and only if $A_i$ is invertible for all $1 \leq i \leq d.$
\end{cor}
\begin{proof}
First, let us assume that $A$ is invertible. Let $B$ be the inverse of $A$ with
\begin{equation*}
B=\left[\begin{array}{c|c|c|c}
B_1 & b_{12}\ones & \cdots & b_{1d}\ones \\
\hline
b_{21}\ones & B_2 & \cdots & b_{2d}\ones \\
\hline
\vdots & \vdots & \ddots & \vdots \\
\hline
b_{d1}\ones & a_{d2}\ones & \cdots & B_d
\end{array}\right].
\end{equation*}
The direct calculation of $AB$ shows that 
\[ A_iB_i = I_{k_i}, \forall \; 1 \leq i \leq d. \]

This shows that $A_i$ is invertible for all $i$. Conversely, suppose that $A_i$ is invertible for all $i$. By scaling $A$ by  a block diagonal matrix, we can assume that $\epsilon(A_i)=1$.  Let $B_i$ be the inverse of $A_i$. We can see that $\epsilon(B_i)=1.$ Equation \ref{eq:multiplicationAB} shows that the following matrix is the inverse of $A.$

\begin{equation*}
B=\left[\begin{array}{c|c|c|c}
B_1 & -a_{12}\ones & \cdots & -a_{1d}\ones \\
\hline
-b_{21}\ones & B_2 & \cdots & -a_{2d}\ones \\
\hline
\vdots & \vdots & \ddots & \vdots \\
\hline
-a_{d1}\ones & -a_{d2}\ones & \cdots & B_d
\end{array}\right].
\end{equation*}
\end{proof}

We end this section with a formula for the number of units in the join algebra for the modular case.  

\begin{cor}
Let $G_i$, for $1 \le i \le d$, be  finite $p$-groups, and let $k$ be a finite field of characteristic $p$. Then we have
\[|(\sJ_{G_1, \ldots, G_d}(k))^\times| = (|k|-1)^d|k|^{(\sum_i |G_i|)+d^2-2d}.\]
\end{cor}

\begin{proof}
By Theorem \ref{units-modular}, we have 
\[ \sJ^{\times} \cong (k^{\times})^d \ltimes k^{d^2-d} \times \prod_{i=1}^d U_1(k[G_i]) .\]
Since each $G_i$ is a $p$-group, an element $u = \sum_{g \in G_i} \alpha_g g$ is in $U_1(k[G_i])$ if and only if $\epsilon(G_i)= \sum_{g \in G} \alpha_g =1$. The number of such elements is $|k|^{|G_i|-1}$ because there are $|G_i|-1$ degrees of freedom for the coefficients. Hence

\begin{eqnarray*}
 |\sJ^{\times}| & =  &  |k^{\times}|^d |k|^{d^2-d}  \prod_{i=1}^d |U_1(k[G_i])| \\
 & = & |k^{\times}|^d  |k|^{d^2-d}  \prod_{i=1}^d|k|^{|G_i|-1} \\
 & = & (|k|-1)^d|k|^{(\sum_i |G_i|)+d^2-2d}.
\end{eqnarray*}
\end{proof}

The unit groups of join algebras are also a source of  $2$-groups, as shown in the following result.

\begin{cor}
Let $G_i$ be  finite $p$-groups. $(\sJ_{G_1, \ldots, G_d}(\F_p))^\times$ ($d > 1$) is a $2$-group if and only if $p=2$.
\end{cor}
\begin{proof}
From the above formula, we have 
\[|(\sJ_{ G_1, \ldots, G_d}(\mathbb{F}_p))^\times | = (p-1)^d p^{(\sum_i |G_i|)+d^2-2d}. \]
This number is a power of $2$ if and only if $p=2$.
\end{proof}

\section{ The Jacobson radical of $\mathcal{J}_{G_1, G_2, \ldots, G_d}(R)$}
In this section, we will investigate the structure of the Jacobson radical of the join algebra $\sJ_{G_1, G_2, \ldots, G_d}(R)$. We will present here two different approaches to this problem. For the first approach, we will utilize the results from the previous section.

Let $k$ be a field of characteristic $p$ (possibly 0) and let $G_{1},\ldots
,G_{d}$ be finite groups of respective orders $k_{1},\ldots ,k_{d}$.  We
identify the group algebra $k[G_{i}]$ with the algebra of $G_{i}$-circulant matrices over $k$.

Write an element $X$ in the join $\mathcal{J}_{G_1, G_2, \ldots, G_d}(k)$ as
\[
X=\left[ 
\begin{array}{llll}
A_{1} & a_{12}\ones_{k_{1},k_{2}} & \cdots  & a_{1d}\ones_{k_{1},k_{d}}
\\ 
a_{21}\ones_{k_{2},k_{1}} & A_{2} & \cdots  & a_{2d}\ones_{k_{2},k_{d}}
\\ 
\vdots  & \vdots  & \ddots  & \vdots  \\ 
a_{d1}\ones_{k_{d},k_{1}} & a_{d2}\ones_{k_{d},k_{2}} & \cdots  & A_{d}%
\end{array}%
\right] 
\]%
where $A_{i}$ is a $G_{i}$-circulant matrix and $\ones_{a,b}$ is an $a\times b$
matrix, all entries of which are $1$. Writing $\text{Rad}(R)$ for the Jacobson
radical of a ring $R$, we then have

\begin{thm}
For $X$ as above, $X\in \Rad(\mathcal{J}_{G_1, G_2, \ldots, G_d}(k))$ if
and only if $A_{i}\in \Rad(k[G_{i}]),i=1,\ldots ,d$ and $a_{ij}=0$
whenever $p\nmid k_{i}k_{j},1\leq i\neq j\leq d$.
\end{thm}

\begin{proof}
Without loss of generality, we reorder $G_{1},\ldots ,G_{d}$ so that $p\nmid
k_{i}$ for $i\leq r$ but $p\mid k_{i}$ for $i>r$. 

First suppose $X\in \text{Rad}(\mathcal{J}_{G_1, G_2, \ldots, G_d}(k))$.  This
implies that $\epsilon(X)\in \text{Rad}(\text{Im}(\epsilon ))$, so we should compute $%
\text{Im}(\epsilon )$.  A little thought reveals that a typical element
of $\text{Im}(\epsilon )$ has the form $ 
\begin{bmatrix}
B & 0 \\ 
C & D%
\end{bmatrix}$ where%
\begin{eqnarray*}
&&B\text{ is an arbitrary }r\times r\text{ matrix} \\
&&C\text{ is an arbitrary }(d-r)\times r\text{ matrix, and} \\
&&D\text{ is an arbitrary diagonal }(d-r)\times (d-r)\text{ matrix}.
\end{eqnarray*}%
The set of all matrices of this form admits a projection homomorphism $\pi $
onto $M_{r}(k)\oplus k^{d-r}$ (taking the above matrix to $(B,%
\vec{v})$ where $\vec{v}$ is the vector of diagonal
entries of $D$), the Jacobson radical of which is $0$.  It follows that $%
X\in \text{Rad}(\mathcal{J}_{G_1, G_2, \ldots, G_d}(k))\Longrightarrow \pi \circ
\epsilon (X)=0$; this immediately gives us that $\epsilon
_{i}(A_{i})=0,1\leq i\leq d$, and $a_{ij}=0$ for $1\leq i\neq j\leq r$%
, i.e. whenever $p\nmid k_{i}k_{j}.$ It remains to see that $A_{i}\in
\text{Rad}(k[G_{i}])$.

We use the characterization $X\in \text{Rad}(\mathcal{J}_{G_1, G_2, \ldots,
G_d}(k))$ if and only if $1+XY$ is a unit in $\mathcal{J}_{G_1, G_2,
\ldots, G_d}(k)$ for all $Y\in \mathcal{J}_{G_1, G_2, \ldots, G_d}(k)$. \
Take $Y$ of the form 
\[
Y=\left[
\begin{array}{ll}
B &  \\ 
0 & \ast  \\ 
\vdots  &  \\ 
0 & 
\end{array}%
\right].
\]%
Then $1+XY$ will have an upper leftmost block $I_{d_{1}}+A_{1}B$, which
according to our characterization of units, must be invertible or almost
invertible (for all $B$).  But this matrix can never be almost invertible
as its augmentation is $1$ (since $\epsilon _{1}(A_{1})$ is now known to
be $0$); thus in fact $A_{1}\in \text{Rad}(k[G_{1}])$ by the same characterization
of Rad.  Similarly we see $A_{i}\in \text{Rad}(k[G_{i}])$ generally. 

Now suppose that the conditions $A_{i}\in \text{Rad}(k[G_{i}]),i=1,\ldots ,d$ and $%
\alpha _{ij}=0$ whenever $p\nmid k_{i}k_{j},1\leq i\neq j\leq d$ hold; we
must show $X\in \text{Rad}(\mathcal{J}_{G_1, G_2, \ldots, G_d}(k))$.  We will
see that $1+XY$ is a unit for any $Y$.  Set $\ $%
\[
Y=\left[ 
\begin{array}{llll}
A_{1}^{\prime } & a_{12}^{\prime }\ones_{k_{1},k_{2}} & \cdots  & a_{1d}^{\prime }\ones_{k_{1},k_{d}} \\ 
a_{21}^{\prime }\ones_{k_{2},k_{1}} & A_{2}^{\prime } & \cdots  & a_{2d}^{\prime }\ones_{k_{2},k_{d}} \\ 
\vdots  & \vdots  & \ddots  & \vdots  \\ 
a_{d1}^{\prime }\ones_{k_{d},k_{1}} & a_{d2}^{\prime
}\ones_{k_{d},k_{2}} & \cdots  & A_{d}^{\prime }%
\end{array}%
\right] 
\]%
The diagonal blocks in $1+XY$ have the form $I+A_{i}A_{i}^{\prime }$ if $%
i\leq r$ or $I+A_{i}A_{i}^{\prime }+wJ_{k_{i},k_{i}}$ if $i>r$.  Since $%
A_{i}\in \text{Rad}(k[G_{i}]),$ $I+A_{i}A_{i}^{\prime }$ is invertible; adding the
(commuting) nilpotent matrix $wJ_{d_{i},d_{i}}\ $if $i>r$ will not change
that.  Thus the diagonal blocks are all invertible. It will follow from
our characterization of units that if, in addition, $\epsilon (1+XY)$ is an
invertible matrix, then $1+XY$ is a unit. We know that $\epsilon (1+XY)$
has the form $\begin{bmatrix}
B & 0 \\ 
C & D%
\end{bmatrix}$ as above, and that all entries on the main diagonal are nonzero. 
(For of course the condition $A_{i}\in \text{Rad}(k[G_{i}])$ forces $\epsilon
_{i}(A_{i})=0$.)  Now consider the $i,j$ block off the diagonal in $1+XY$
where $i,j\leq r$.  We will see that it is actually a zero matrix.  This
block is a sum of terms (i) $a_{il}a_{lj}^{\prime }k_{l}J_{ij}\
(i\neq l\neq j)$, (ii) $\epsilon _{i}(A_{ii})a_{ij}^{\prime }J_{ij}$%
, and (iii) $a_{ij}\epsilon _{j}(A_{jj}^{\prime })J_{ij}$.  Terms
of type (i) are all zero since\ either $\alpha _{il}=0\ (l\leq r)$ or $%
k_{l}=0\ (l>r)$, the term (ii) is zero since $\epsilon _{i}(A_{ii})=0$,
and finally term (iii) is zero since $a_{ij}=0$.  Thus we see, taking 
$\epsilon (1+XY)$, that the matrix $B$ is diagonal, and the entire matrix 
$\epsilon (1+XY)$ is, in particular, lower triangular with nonzero diagonal
entries. Thus it is invertible, and the criteria for invertibility of $1+XY
$ are met.
\end{proof}

\begin{cor}
Let $G_i$ ($1 \le i \le d$)   be finite $p$-groups and let $k$ be a finite field of characteristic $p$. Then we have 
\[ | \Rad(\sJ_{G_1, \ldots, G_d}(k) | = |k|^{\Sigma |G_i|+d^2-2d}. \]
\end{cor}

\begin{proof}
An element 
\[
X=\left[ 
\begin{array}{llll}
A_{1} & a_{12}\ones_{k_{1},k_{2}} & \cdots  & a_{1d}\ones_{k_{1},k_{d}}
\\ 
a_{21}\ones_{k_{2},k_{1}} & A_{2} & \cdots  & a_{2d}\ones_{k_{2},k_{d}}
\\ 
\vdots  & \vdots  & \ddots  & \vdots  \\ 
a_{d1}\ones_{k_{d},k_{1}} & a_{d2}\ones_{k_{d},k_{2}} & \cdots  & A_{d}%
\end{array}
\right] 
\] 
in $\sJ_{G_1, \ldots, G_d}(k)$ belongs to the Jacobson radical if and only if $\epsilon(A_i) =  0$ for all $i$  and with no restriction on  the off-diagonal blocks. Clearly, $\epsilon(A_i) = 0$ if and only if the row sum of $A_i$ is zero. This means we have $|G_i|-1$ degrees of freedom which gives $|k|^{|G_i|-1}$ elements in $A_i$.  Since there are no  restrictions on the off-diagonal blocks, we have a total of 
\[|k|^{|G_1|-1}\ldots |k|^{|G_d|-1} |k|^{d^2-d} = |k|^{\Sigma |G_i|+d^2-2d}\]
elements in the Jacobson radical.
\end{proof}

\begin{cor} $\sJ_{G_1, G_2, \ldots, G_d}(k)$ is semisimple if and only if $|G_i|$ are invertible in $k$. 
\end{cor}
We discuss another approach to this problem, which may be of independent interest. Let $G_1, G_2, \ldots, G_d$ be as before. We will work with a general ring $R$ that satisfies the following Hypothesis. 

\begin{Hypothesis} $R$ is a semisimple ring  and for each $1 \leq i \leq d$, either $|G_i|=0$ in $R$ or $|G_i|$ is invertible. 
\end{Hypothesis}
In particular, a field would automatically satisfy this condition.  If $|G_i|$ is invertible in $R$ for all $1 \leq i \leq d$ then by Theorem \ref{thm:maschke},  $\sJ_{G_1, G_2, \ldots, G_d}(R)$ is semisimple so its Jacobson radical is $0$. Therefore, we can assume that, up to order, there exists a (unique) positive integer $r$ such that 

\begin{itemize}
    \item $|G_i|$ is invertible in $R$ for $1 \leq i \leq r$. 
    \item $|G_i|=0$ in $R$ for $r < i \leq d$.
\end{itemize}
Let us first explain our strategy for this second approach. We will find an ideal $\Delta$ - as small as possible - such that the quotient ring $\sJ_{G_1, G_2, \ldots, G_d}(R)/\Delta$ is semisimple. This, in turn, can be done by constructing a surjective ring homomorphism from $\sJ_{G_1, G_2, \ldots, G_d}(R)$ to another semisimple ring. Our strategy is based on the following observation. 
\begin{prop} (\cite[Section 4.3, Lemma b]{pierce1982associative}) \label{prop:rad_quotient}
Let $R$ be a ring and $I$ a two-sided ideal of $R$ such that $R/I$ is semisimple. Let $\Rad(R)$ be the Jacobson radical of $R$. Then $\Rad(R) \subseteq I.$
\end{prop}
First, we discuss an elementary lemma. 
\begin{lem}
Let $R$ be a semisimple ring, $G$ a group with either $|G|=0$ in $R$ or $|G|$ invertible in $R$ and let
\[ e_{G}=\sum_{g \in G} g.\] 
Then $e_G$ is an element of the Jacobson radical of $R[G]$ if and only if $|G|=0$ in $R.$
\end{lem}
\begin{proof}
If $|G|$ is invertible in $R$, then by Maschke's theorem, $R[G]$ is semisimple, so the Jacobson radical of $R[G]$ is $0$. Therefore, $e_G$ cannot be an element of the Jacobson radical of $R[G].$ Conversely, assume that $|G|=0$ in $R.$ We claim that $1+e_G y$ is a unit in $R[G]$ for all $y \in R[G].$  We have 
\[ e_G y =\epsilon(y) e_G, \]
where $\epsilon(y)$ is the augmentation of $y$. In particular, we have 
\[ e_G^2= |G| e_G=0 .\] 
Therefore $(e_Gy)^2=\epsilon(y)^2 e_G^2=0.$ This shows that $1+e_G y$ is invertible. In fact, its inverse is exactly $1-e_Gy.$
This shows that $e_G$ belongs to the Jacobson radical of $R[G].$
\end{proof}
For each $1 \leq i \leq d$, let $I_i$ be the Jacobson radical of the group ring $R[G_i].$ Note that by Maschke's theorem \ref{thm:classical_maschke}, for $1 \leq i \leq r$, $R[G_i]$ is semisimple, so $I_i=0$. Also from our assumption that $R$ is semisimple, $R[G_i]/I_i$ is semisimple for all $1 \le i \le d$.

Let us consider a generic element of $\sJ_{G_1, G_2, \ldots, G_d}(R)$ 

\begin{equation*}
A=\left[\begin{array}{c|c|c|c}
C_1 & a_{12}\ones & \cdots & a_{1d}\ones \\
\hline
a_{21}\ones & C_2 & \cdots & a_{2d}\ones \\
\hline
\vdots & \vdots & \ddots & \vdots \\
\hline
a_{d1}\ones & a_{d2}\ones & \cdots & C_d
\end{array}\right].
\end{equation*}
We can further partition $A$ into the following blocks 
\[ A =\begin{bmatrix} A_1 & B_1 \\ B_2 & A_2 \end{bmatrix},\] 
where $A_1$ is the union of the upper $r$ blocks, $A_2$ is the union of the lower $d-r$ blocks, $B_1$ (respectively $B_2$) is the union of the upper right (respectively lower left) blocks. Concretely, we have 
\begin{equation*}
A_1=\left[\begin{array}{c|c|c|c}
C_1 & a_{12}\ones & \cdots & a_{1r}\ones \\
\hline
a_{21}\ones & C_2 & \cdots & a_{2r}\ones \\
\hline
\vdots & \vdots & \ddots & \vdots \\
\hline
a_{r1}\ones & a_{r2}\ones & \cdots & C_r
\end{array}\right],
\end{equation*}
\begin{equation*}
A_2=\left[\begin{array}{c|c|c|c}
C_{r+1} & a_{r+1,r+2}\ones & \cdots & a_{r+1,d}\ones \\
\hline
a_{r+2,r+1}\ones & C_2 & \cdots & a_{r+2,d}\ones \\
\hline
\vdots & \vdots & \ddots & \vdots \\
\hline
a_{d,r+1}\ones & a_{d,r+2}\ones & \cdots & C_d
\end{array}\right].
\end{equation*}
Similarly for $B_1, B_2$. Note that we can consider $A_1$ (respectively $A_2$) as an element of $\sJ_{G_1, \ldots, G_r}(R)$ (respectively $\sJ_{G_{r+1}, \ldots, G_d}(R)$.) Suppose $X$ is another element in $\sJ_{G_1, \ldots, G_d}(R)$ of the form 

\begin{equation*}
X=\left[\begin{array}{c|c|c|c}
D_1 & x_{12}\ones & \cdots & x_{1d}\ones \\
\hline
x_{21}\ones & D_2 & \cdots & x_{2d}\ones \\
\hline
\vdots & \vdots & \ddots & \vdots \\
\hline
x_{d1}\ones & x_{d2}\ones & \cdots & D_d
\end{array}\right].
\end{equation*}
Again, we can write $X$ in the following form

\[ X =\begin{bmatrix} X_1 & Y_1 \\ Y_2 & X_2 \end{bmatrix} ,\]
then we have 
\[ AX= \begin{bmatrix} A_1X_1+B_1 Y_2 & A_1 Y_1 + B_1 X_2 \\ B_2 X_1+A_2 Y_2 & B_2 Y_1+ A_2 X_2 \end{bmatrix} .\] 
We note that $B_1 Y_2=0$ and $B_2 Y_1$ consists of the blocks of form $c J_{m,n}$ for suitable $m$, $n$, and $c$. We also note that the diagonal blocks of $A_2X_2$ are just the $C_iD_i$'s, $r+1\le i \le d$.

\begin{prop}
Let $\psi$ be the map 
\[ \psi\colon \sJ_{G_1, G_2, \ldots, G_d}(R) \to \sJ_{G_1, \ldots, G_r}(R) \times \prod_{r+1 \leq i \leq d} R[G_i]/I_i ,\]
sending 

\[ A \mapsto (A_1, \overline{C_{r+1}}, \ldots, \overline{C_{d}}) .\] 
Then $\psi$ is a surjective ring homomorphism. 
\end{prop}

\begin{proof}
Let $A,X \in \sJ_{G_1, \ldots, G_d}(R)$ as before. We need to show that
\[ \psi(A+X)=\psi(A)+\psi(X), \]
and 
\[ \psi(AX)=\psi(A) \psi(X).\] 
The first identity is obvious. Let us focus on the second identity. By the above calculations and the fact that $e_{|G_i|} \in I_i$ for $r+1 \leq i \leq d$, we see that 
\[ \psi(AX)=(A_1X_1, \overline{C_{r+1} D_{r+1}}, \ldots, \overline{C_{d} D_d}) =\psi(A) \psi(X) .\] 
We conclude that $\psi$ is a ring homomorphism. It is surjective because for an element \\ $(A_1, \overline{C_{r+1}}, \ldots, \overline{C_{d}}) \in \sJ_{G_1, \ldots, G_r}(R) \times \prod_{r+1 \leq i \leq d} R[G_i]/I_i$, we have  
\[ \psi(A)= (A_1, \overline{C_{r+1}}, \ldots, \overline{C_{d}}), \]
where 
\[ A= A_1 \oplus C_{r+1} \oplus \ldots \oplus C_d\]
is a diagonal block matrix whose block components are $A_1, C_{r+1}, \ldots, C_d$.
\end{proof}
Let $\Delta$ be the kernel of this ring homomorphism. Then we have 
\[ \sJ_{G_1, \ldots, G_d}(R)/ \Delta \cong \sJ_{G_1, \ldots, G_r}(R) \times \prod_{i=r+1}^d R[G_i]/I_i ,\] 
By the generalized Maschke's theorem \ref{thm:maschke}, $\sJ_{G_1, \ldots, G_d}(R) \times \prod_{i=r+1}^d R[G_i]/I_i$ is a semisimple ring, we conclude that  $I \subseteq \Delta$, where $I$ is the Jacobson radical of $\sJ_{G_1, G_2, \ldots, G_d}(R).$
  We will now prove the other direction, namely $\Delta \subseteq I.$ To do so, we use the following lemma. 
\begin{lem} (see \cite[Theorem 4.12]{lam}) \label{lem:nil}
Let $A$ be a left Artinian algebra. The $\Rad(A)$ is the largest nilpotent left ideal and it is also the largest nilpotent right ideal. In particular 
\begin{enumerate}
    \item $\Rad(A)$ is a two-sided nilpotent ideal. 
    \item If $M$ is a two-sided nilpotent ideal then $M \subseteq \Rad(A).$
\end{enumerate}
\end{lem}

We will use the second statement to show that $\Delta \subseteq I$. Namely, we will show that all elements of $\Delta$ are nilpotent. Let $A \in \Delta$. Then as before, $A$ has the following form 
\[ A =\begin{bmatrix} 0 & B_1 \\ B_2 & A_2 \end{bmatrix}.\] 
where $B_1, B_2, A_2$ are as before. Furthermore, if we write $A_2$ in the form 
\begin{equation*}
A_2=\left[\begin{array}{c|c|c|c}
C_{r+1} & a_{r+1,r+2}\ones & \cdots & a_{r+1,d}\ones \\
\hline
a_{r+2,r+1}\ones & C_{r+2} & \cdots & a_{r+2,d}\ones \\
\hline
\vdots & \vdots & \ddots & \vdots \\
\hline
a_{d,r+1}\ones & a_{d,r+2}\ones & \cdots & C_d
\end{array}\right],
\end{equation*}
then we must have $C_i \in I_i$ for $r+1 \leq i \leq d.$ In particular, $C_i$ are all nilpotent (by Lemma \ref{lem:nil}). We also note that $\epsilon(C_i)=0.$ In fact, by Proposition \ref{prop:rad_quotient}, $\Rad(R[G_i]) \subseteq \ker(\epsilon)$ since $R[G_i]/\ker(\epsilon) \cong R$ which is semisimple by our assumption. Direct calculations show that $B_1A_2=A_2B_2=0$ and hence $A^2$ is of the following form 
\[ A^2 =\begin{bmatrix} 0 & 0 \\ 0 & A'_2 \end{bmatrix},\]
where 
\[ A'_2=C_{r+1}^2 \oplus \ldots \oplus C_{d}^2 .\] 
Since $C_i$ are all nilpotent, we conclude that $A^2$ is nilpotent, and hence $A$ is nilpotent as well. This shows that $\Delta \subseteq I.$ 

In summary, we have 
\begin{thm} \label{thm:jacobson}
The Jacobson radical of $\sJ_{G_1, G_2, \ldots, G_d}(R)$ is the kernel of the surjective ring homomorphism 
\[ \psi\colon \sJ_{G_1, G_2, \ldots, G_d}(R) \to \sJ_{G_1, \ldots, G_r}(R) \times \prod_{r+1 \leq i \leq d} R[G_i]/I_i .\]
Concretely, let 

\begin{equation*}
A=\left[\begin{array}{c|c|c|c}
C_1 & a_{12}\ones & \cdots & a_{1d}\ones \\
\hline
a_{21}\ones & C_2 & \cdots & a_{2d}\ones \\
\hline
\vdots & \vdots & \ddots & \vdots \\
\hline
a_{d1}\ones & a_{d2}\ones & \cdots & C_d
\end{array}\right].
\end{equation*}
Then $A$ belongs to the Jacobson radical of $\sJ_{G_1, G_2, \ldots, G_d}(R)$ if and only if the following conditions are satisfied: 
\begin{enumerate}
    \item $C_i=0, 1 \leq i \leq r$, 
    \item $a_{ij}=0, 1 \leq i,j \leq r$, 
    \item $C_i \in I_i, r+1 \leq i \leq d$.
\end{enumerate}
\end{thm}

\begin{cor}
Suppose that $R=k$ is an algebraically closed field of characteristics $p$. Let $G_1, G_2, \ldots, G_d$ be as before. Then, the number of irreducible modules over $\sJ_{G_1, G_2, \ldots, G_d}(k)$ is 
\[ \sum_{i=1}^{d} c_p(G)- r +1,\] 
where $c_p(G_i)$ is the number of $p$-regular conjugacy classes of $G_i$.
\end{cor}
\begin{proof}
For a ring $R$, we define the semisimplification of $R$ as 
\[ R^{\se} = R/\Rad(R) .\] 
A simple module over $R$ is of the form $R/\mathfrak{m}$ where $\mathfrak{m}$ is a left maximal ideal in $R$. Since $\Rad(R)$ is the intersection of all left maximal ideals in $R$, we conclude that there is a bijection between the set of simple modules over $R$ and the set of simple modules over $R^{\se}.$  From this observation and the isomorphism discussed in Theorem \ref{thm:jacobson}
we conclude that the number of irreducible modules over $\sJ_{G_1, G_2, \ldots, G_d}(k)$ is the same as the number of irreducible modules over $\sJ_{G_1, G_2, \ldots, G_r}(k) \times \prod_{i=r+1}^d k[G_i]/\Rad(k[G_i])$. By \ref{prop:number_of_irr_modules}, we know that the number of irreducible modules over $\sJ_{G_1, G_2, \ldots, G_r}(k)$ is exactly 
\[ \sum_{i=1}^r c(G_i) -r+1 = \sum_{i=1}^r c_p(G_i)-r+1 .\] 
Additionally, the number of irreducible modules over $k[G_i]/\Rad(k[G_i])$ is the same as the number of irreducible modules over $k[G_i]$ which is known to be $c_p(G_i)$ (see \cite{reiner1964number}). We conclude that the number of irreducible modules over $\sJ_{G_1, G_2, \ldots, G_d}(k)$ is exactly 
\[ \sum_{i=1}^d c_p(G_i)-r+1.
\qedhere\]

\end{proof}

\section{The join algebra $\sJ_{G_1, G_2, \ldots, G_d}(k)$ and Frobenius algebras}

An important class of algebras is Frobenius algebra which we now recall. 
\begin{defn} (\cite[Page 66-67]{lam2012lectures})
Let $A$ be a finite-dimensional $k$-algebra. Then $A$ is called a Frobenius algebra if one of the following equivalent conditions holds 
\begin{enumerate}
    \item There exists a non-degenerate bilinear form $\sigma: A \times A \to k$ such that $\sigma(ab,c)=\sigma(a,bc)$ for all $a,b, c \in A.$ Here non-degenerate means that if $\sigma(x,y)=0$ for all $x$ then $y=0.$ We call $\sigma$ a Frobenius form of $A.$
    \item There exists a linear map $\lambda: A \to k$ such that the kernel of $\lambda$ contains no nonzero left ideal of $A.$
\end{enumerate}
\end{defn}

It is known that if $k$ is a field and $G$ is a finite group, then the group algebra $k[G]$ is always a Frobenius algebra regardless of the characteristic of the field $k$ (see \cite[Example 3.15E]{lam2012lectures}). In this section, we completely answer the following question: when is the join algebra $\sJ_{G_1, G_2, \ldots, G_d}(k)$ a Frobenius algebra? 
\begin{thm} \label{thm:frob}
Suppose $G_1, G_2, \ldots, G_d$ are groups over a field $k$ of characteristic $p$ with $d \geq 2$ Then, the join algebra $\sJ_{G_1, G_2, \ldots, G_d}(k)$ is a Frobenius algebra if and only if $|G_i|$ is invertible in $k$ for all $1 \leq i \leq d.$ 
\end{thm}
We remark that if $|G_i|$ are invertible in $k$ then by Theorem \ref{thm:decomposition}, $\sJ_{G_1, G_2, \ldots, G_d}(k)$ is semisimple, hence a Frobenius algebra by \cite[Example 3.15D]{lam2012lectures}.  Therefore, it is sufficient to consider the case that at least one of  $|G_i|$ is $0$ in $k$. Without loss of generality, we can assume that $|G_1|=0$ in $k.$ Our key observation is that there are many left ideals in $\sJ_{G_1, G_2, \ldots, G_d}(k).$

Let $(a_1,a_2, \ldots, a_d) \in k^d$. We define
\[ v_{a_1, a_2, \ldots, a_d}= 
\left(\begin{array}{c|c|c|c}
a_1 J_{k_1, k_1} & a_{2} J_{k_1, k_2} & \cdots & a_{d}J_{k_1, k_d} \\
\hline
0  & 0 & \cdots & 0 \\
\hline
\vdots & \vdots & \ddots & \vdots \\
\hline
0 & 0  & \cdots &0 
\end{array}\right) .\] 
 
Let $I_{a_1,a_2, \ldots, a_d}$ be the vector space generated by $v_{a_1,a_2, a_3, \ldots, a_d}$. Then 
\begin{prop} \label{prop:7.3}
$I_{a_1,a_2, \ldots, a_d}$ is a left ideal of $\sJ_{G_1, G_2, \ldots, G_d}(k)$. If $(a_1,a_2, \ldots, a_d) \neq (0,0, \ldots, 0)$ then $I_{a_1,a_2, \ldots, a_d}$ is not $0$.
\end{prop}
\begin{proof}
Let 
\[ A=\left(\begin{array}{c|c|c|c}
C_1 & a_{12}\ones & \cdots & a_{1d}\ones \\
\hline
a_{21}\ones & C_2 & \cdots & a_{2d}\ones \\
\hline
\vdots & \vdots & \ddots & \vdots \\
\hline
a_{d1}\ones & a_{d2}\ones & \cdots & C_d
\end{array}\right) .\] 
Then 
\[ A v_{a_1, a_2, \ldots, a_d} = \epsilon(C_1) v_{a_1, a_2, \ldots, a_d} \in I_{a_1, a_2, \ldots, a_d}.\] 
\end{proof}

\begin{prop} \label{prop:7:4}
Let $\lambda: \sJ_{G_1, G_2, \ldots, G_d}(k) \to k$ be a linear functional. Then there exists $(a_1, \ldots, a_d) \neq (0, 0, \ldots, 0)$ such that $\lambda(v_{a_1, a_2, \ldots, a_d})=0.$ Consequently $\lambda(I_{a_1, a_2, \ldots, a_d})=0.$
\end{prop}

\begin{proof}
Let $V$ be the $d$-dimensional vector space 
\[ V = \left\{ \left(\begin{array}{c|c|c|c}
x_1 J_{k_1, k_1} & x_{2} J_{k_1, k_2} & \cdots & x_{d}J_{k_1, k_d} \\
\hline
0  & 0 & \cdots & 0 \\
\hline
\vdots & \vdots & \ddots & \vdots \\
\hline
0 & 0  & \cdots &0 
\end{array}\right)| (x_1, x_2, \ldots, x_d) \in k^d  \right \} .\] 

The restriction of $\lambda$ to $V$ induces a linear functional map 
\[ \lambda: V \to k .\] 
If $d>1$, this map must have a non-trivial kernel. So, there must exist $(a_1, a_2, \ldots, a_d) \neq 0$ such that $\lambda(v_{a_1, a_2, \ldots, a_d})=0.$
\end{proof}
We see that Theorem \ref{thm:frob} is then a consequence Proposition \ref{prop:7.3} and Proposition \ref{prop:7:4}. Here are two direct corollaries of this theorem. 
\begin{cor}
Assume that $d \geq 2$ and $p| |G_1|$. Then $\sJ_{G_1, G_2, \ldots, G_d}(k)$ cannot be the group algebra of any finite group $G.$
\end{cor}

In Example \ref{ex:join}, we discuss the possibility of writing a group algebra as a join algebra. It turns out that this is not possible in the modular case. 
\begin{cor}
Let $G$ be a group such that $|G|=0$ in k. Then $k[G]$ is not isomorphic to any $\sJ_{G_1, G_2, \ldots, G_d}(k)$ with $d \geq 2.$
\end{cor}
\begin{proof}
Assume to the contrary that $k[G]=\sJ_{G_1, G_2, \ldots, G_d}(k)$ for some $d \geq 2.$ Since $k[G]$ is a Frobenius algebra, $\sJ_{G_1, G_2, \ldots, G_d}(k)$ is a Frobenius algebra as well. By the above corollary, $|G_i|$ must be invertible in $k$ for $1 \leq i \leq d.$ By Theorem \ref{thm:maschke}, this implies that $\sJ_{G_1, G_2, \ldots, G_d}(k)$ is semisimple. However, since $|G|=0$ in $k$, $k[G]$ is not semisimple. This leads to a contradiction. 

\end{proof}

\section{Artin-Wedderburn decomposition/Generalized Circulant Diagonalization Theorem}
In this section, we describe an explicit isomorphism mentioned in Example \ref{ex:cyclic}. Throughout this section, we will assume that $G_i = \Z/k_i$ is a cyclic group of order $k_i$. In this section, we will assume that $R=k$ is a field (for applications that we have in mind, $k=\C$ would be sufficient). Additionally, for simplicity, we will use the notation  $\sJ_{k_1, k_2, \ldots, k_d}(k)$ for $\sJ_{G_1, G_2, \ldots, G_d}(k).$

We recall that the $n$-dimensional Discrete Fourier Transform (DFT) is the linear map $\C^n\to \C^n$ represented in matrix form by the \textit{DFT matrix}
\[
F_n= \begin{pmatrix}
1&1&1&\cdots &1 \\
1&\omega&\omega^2&\cdots&\omega^{n-1} \\
1&\omega^2&\omega^4&\cdots&\omega^{2(n-1)}\\
\vdots&\vdots&\vdots&\ddots&\vdots\\
1&\omega^{n-1}&\omega^{2(n-1)}&\cdots&\omega^{(n-1)(n-1)}
\end{pmatrix},
\]
where $\omega$ is a primitive $n$-th root of unity. The matrix $F_n$ is invertible, with inverse $F_n^{-1}=\frac{1}{n}F_n^*$ (where $M^*$ denotes the conjugate transpose of the matrix $M$).
Moreover, the Circulant Diagonalization Theorem states that all the circulant matrices of size $n$ can be simultaneously diagonalized by conjugation with $F_n$ (see \cite[Theorem 3.2.1]{davis2013circulant}).
Note that, although the DFT matrix and the Circulant Diagonalization Theorem for $n \times n$ circulant matrices are usually introduced over the field of complex numbers, they make sense over any field $k$ containing a primitive $n$-th root of unity and in which $n$ is invertible. 

The Generalized Circulant Diagonalization Theorem of \cite{CM1} can be thought of as saying that, for any $d,k_1,\dots,k_d\in\N$, it is possible to turn all matrices of $\sJ_{k_1,\dots,k_d}(\C)$ into an almost diagonal form by conjugation with a block-diagonal matrix whose diagonal blocks are DFT matrices of suitable size.
\begin{defn}
For $d,k_1,\dots,k_d\in\N$, the \textit{join-DFT} matrix of sizes $k_1,\dots,k_d$ is the block-diagonal matrix
\[
F_{k_1,\dots,k_d}=\left(\begin{array}{c|c|c|c}
F_{k_1} & \mathbf{0} & \dots & \mathbf{0}\\
\hline
\mathbf{0} & F_{k_2} & \dots & \mathbf{0}\\
\hline
\vdots & \vdots & \ddots & \vdots\\
\hline
\mathbf{0} & \mathbf{0} & \dots & F_{k_d}
\end{array}
\right).
\]
To keep the notation cleaner, we will also use the shorthand $F_{\mathbf{k}}=F_{k_1,\dots,k_d}$. We recall that for any $n  \in \mathbb{N}$   and a field $k$, a primitive  $n$th of unity is an element in the multiplicative group $k^{\times}$ of $k$, that generates 
a subgroup of $k^{\times}$ containing all distinct roots of the polynomial $x^n - 1$. This definition is less restrictive than asking this group to contain $n$ elements. 
In particular, the existence of a $n$th-primitive root does not imply that $n$ is invertible in $k$.
\end{defn}

\begin{thm}\label{thm:Artin-Wedderburn}
Let $d,k_1,\dots,k_d\in\N$. If $k$ contains the inverses of $k_1,\dots,k_d$ and the primitive roots of unity of orders $k_1,\dots,k_d$, then the algebra $\sJ_{k_1,\dots,k_d}(k)$ is isomorphic to \[\underbrace{k\times\dots\times k}_{\text{$k_1+\dots+k_d-d$}}\times\ M_d(k).\] The isomorphism is given by conjugation with the product of the join-DFT matrix $F_{k_1,\dots,k_d}$ and a suitable permutation matrix.
\end{thm}
\begin{proof}
The matrix $F_{\mathbf{k}}=F_{k_1,\dots,k_d}$ contains exactly $d$ columns whose nonzero entries are all equal. With our conventions, these columns are the first, the $(k_1+1)$-th, the $(k_1+k_2+1)$-th, ..., the $(k_1+\dots+k_{d-1}+1)$-th, that is, the columns containing the first column of each diagonal block. Let us refer to these as the \textit{bad} columns, and to the others as the \textit{good} ones.
By \cite[Theorem 1]{CM1}, the good columns of $F_\mathbf{k}$ are common eigenvectors of all matrices in $\sJ_{k_1,\dots,k_d}(R)$.
Let $P$ be the permutation matrix that brings the bad columns at the end of $F_\mathbf{k}$, otherwise keeping the relative order of both the good and bad columns. Then, for all $A\in \sJ_{k_1,\dots,k_d}(k)$ of the form \ref{eq:join circulant matrix}, the matrix $P^{-1}F_\mathbf{k}^{-1}AF_\mathbf{k}P$ has the shape
\[D_A\oplus \overline{A}=
\left(\begin{array}{c|c}
    D_A & \mathbf{0} \\
    \hline
    \mathbf{0}  & \overline{A}
\end{array}\right),
\]
where $D_A$ is the diagonal matrix having the \textit{circulant} eigenvalues of $A$ (in the sense of \cite[Definition 1]{CM1}) on the diagonal, and
\[
\overline{A}=
\begin{pmatrix}
\epsilon(C_1) & k_2a_{12} & \dots & k_da_{1d}\\
k_1a_{21} & \epsilon(C_2) & \dots & k_da_{2d}\\
\vdots & \vdots & \ddots & \vdots\\
k_1a_{d1} & k_2a_{d2} & \dots & \epsilon(C_d)
\end{pmatrix}.
\]

Let us define $n=k_1+\dots+k_d$ and consider the $k$-algebra homomorphism 
\begin{equation}\label{eq:Artin-Wedderburn decomp}
    \begin{array}{c}
    \Phi:\sJ_{k_1,\dots,k_d}(k)\to k^{n-d}\times\ M_d(k),\\\\
    \Phi: A\mapsto P^{-1}F_{\mathbf{k}}\,\! ^{-1}AF_\mathbf{k}P=D_A\oplus\overline{A}\mapsto \left((D_A)_{11},(D_A)_{22},\dots,(D_A)_{n-d,n-d},\overline{A}\right).
    \end{array}
\end{equation}

The injectivity of $\Phi$ follows from three properties: the invertibility of $k_1,\dots,k_d$, forcing the implication 
$k_ia_{ji}=0\Rightarrow a_{ij}=0$; the fact that diagonal entries of $D_A\oplus\overline{A}$ are precisely the eigenvalues of the circulant blocks of $A$; and the fact that the circulant blocks of $A$ are diagonalizable, thanks to the presence of the necessary roots of unity in $k$.

As regards the surjectivity of $\Phi$, for all $(r_1,\dots,r_{n-d},M)\in k^{n-d}\times M_d(k)$, the matrices
\[\begin{array}{l}
     C_1=F_\mathbf{k}\cdot diag(M_{11},r_{1},\dots,r_{k_1-1})\cdot F_\mathbf{k}^{-1} \\
     C_2=F_\mathbf{k}\cdot diag(M_{22},r_{k_1},\dots,r_{k_1+k_2-2})\cdot F_\mathbf{k}^{-1} \\ 
     \vdots \\
     C_d=F_\mathbf{k}\cdot diag(M_{dd},r_{k_1+\dots+k_{d-1}-d-2},\dots,r_{n-d})\cdot F_\mathbf{k}^{-1}
\end{array}
\]
are circulant and can be assembled into the join
\[
\left(
\begin{array}{c|c|c|c}
   C_1 & (M_{12}/k_2)\ones & \dots & (M_{1d}/k_d)\ones \\
   \hline
   (M_{21}/k_1)\ones & C_2 & \dots & (M_{2d}/k_d)\ones \\
   \hline
   \vdots & \vdots & \ddots & \vdots \\
   \hline
   (M_{d1}/k_1)\ones & (M_{d2}/k_2)\ones & \dots & C_d
\end{array}
\right),
\]
which is a preimage of $(r_1,\dots,r_{n-d},M)$.
\end{proof}
\begin{rem}
In the literature, there are different conventions for the definition of DFT matrices. Many authors prefer to define the DFT matrix as 
\[\widetilde{F}_n=\frac{1}{\sqrt{n}}F_n.\]
The normalization factor $1/\sqrt{n}$ has the merit of making the DFT a unitary operator. The convention used here has the advantage of requiring less strict assumptions on the field $k$ in Theorem \ref{thm:Artin-Wedderburn}, namely, $k$ need not contain the square roots of $k_1,\dots,k_d$ (and their inverses). We note that the use of different forms of the DFT matrices, provided they are defined over $k$, does not substantially modify the theorem, as the only effect would be to substitute the matrix $\overline{A}$ with a similar matrix.

However, if $k$ happens to contain the square roots of $k_1,\dots,k_d$, the adoption of $\widetilde{F}_n$ as DFT matrix, and the corresponding choice of
\[
\widetilde{F}_{k_1,\dots,k_d}=\left(\begin{array}{c|c|c|c}
\widetilde{F_{k}}_1 & \mathbf{0} & \dots & \mathbf{0}\\
\hline
\mathbf{0} & \widetilde{F_{k}}_2 & \dots & \mathbf{0}\\
\hline
\vdots & \vdots & \ddots & \vdots\\
\hline
\mathbf{0} & \mathbf{0} & \dots & \widetilde{F_{k}}_d
\end{array}
\right)
\]
as a join-DFT matrix have an interesting graph-theoretic consequence. In fact, the adjacency matrix of a graph join of circulant (unweighted) graphs is a join of circulant matrices $A$ as in \eqref{eq:join circulant matrix}, but with all $a_{ij}=1$. Now, the conjugation with $\widetilde{F}_{k_1,\dots,k_d}P$ produces the block-diagonal matrix $D_A\oplus \widetilde{A}$, with $D_A$ as in the theorem, but
\[
\widetilde{A}=\begin{pmatrix}
\epsilon(C_1) & \sqrt{k_1k_2} & \dots & \sqrt{k_1k_d}\\
\sqrt{k_2k_1} & \epsilon(C_2) & \dots & \sqrt{k_2k_d}\\
\vdots & \vdots & \ddots & \vdots\\
\sqrt{k_dk_1} & \sqrt{k_dk_2} & \dots & \epsilon(C_d)
\end{pmatrix}
\]
 is symmetric. Consequently, the adjacency matrix of any graph join of circulant (unweighted) graphs is diagonalizable. Of course, the same is true for more general matrices $A$ as in \eqref{eq:join circulant matrix} having $a_{ij}=a_{ji}$.
\end{rem}

\begin{cor}
In the same hypotheses of Theorem \ref{thm:Artin-Wedderburn}, the map $\overline{\Phi}:\sJ_{k_1,\dots,k_d}(k)\to M_d(k)$, $\overline{\Phi}:A\mapsto\overline{A}$ is an $k$-algebra epimorphism.
\end{cor}

\bibliographystyle{abbrv}
\bibliography{CM2.bib}
\end{document}